\documentclass[10pt]{article}\usepackage{ amsmath, amsthm, amssymb, amsfonts}
\usepackage{pstricks}
\usepackage{pst-all}
\usepackage{multimedia}
\usepackage{dsfont}
\usepackage{comment}
\usepackage{cleveref}
 \usepackage{graphicx}
 \usepackage{caption}
\usepackage{subcaption}
\newcommand\inv{^{-1}}

\renewcommand\a{\alpha}
\renewcommand\b{\beta}
 \newcommand\dl{\delta}
 \newcommand\D{\Delta}
\newcommand{\s}{\sigma} 

\newcommand{\e}{\eta}
\renewcommand\t{\tau}

\newcommand\ph{\varphi}
\newcommand\ve{\varepsilon}

\newcommand{\Om}{\Omega}
 
\newcommand\cd{\cdots}
\newcommand\ld{\ldots}
 \renewcommand\l{\ell}
\newcommand\pl{\partial}
\newcommand\iy{\infty}
\newcommand{\x}{\xi}
\newcommand{\tx}{\tilde{\xi}}

\newcommand{\z}{\zeta}
\newcommand{\tz}{\tilde{\zeta}}
 \def\({\left(} \def\){\right)}

\newcommand{\ra}{\rightarrow}

\newcommand{\bP}{\mathbb{P}}

\newcommand{\bZ}{\mathbb{Z}}
\newcommand{\bC}{\mathbb{C}}

\newcommand{\bN}{\mathbb{N}}

\newcommand{\cB}{\mathcal{B}}

\newcommand{\cC}{\mathcal{C}}
\newcommand{\cX}{\mathcal{X}}
\newcommand{\cZ}{\mathcal{Z}}

\newcommand{\cS}{\mathcal{S}}
\newcommand{\cH}{\mathcal{H}}

\newcommand{\cP}{\mathcal{P}}

\newcommand{\cL}{\mathcal{L}}
\newcommand{\cF}{\mathcal{F}}
\newcommand{\cN}{\mathcal{N}}

\newcommand{\me}{\textrm{e}}
\newcommand{\mi}{\textrm{i}}

\newcommand{\ua}{\uparrow}
\newcommand{\da}{\downarrow}
\newcommand{\hf}{\frac{1}{2}}

\newcommand{\be}{\begin{equation}}
\newcommand{\ee}{\end{equation}}
\newcommand{\ba}{\begin{eqnarray*}}
\newcommand{\ea}{\end{eqnarray*}}
\newcommand{\bae}{\begin{eqnarray}}
\newcommand{\eae}{\end{eqnarray}}
\newcommand{\bc}{\begin{center}}
\newcommand{\ec}{\end{center}}

\newcommand{\fr}{\frac}

\newcommand{\sq}{\sqrt}
\newcommand{\sgn}{\textrm{sgn}}

\begin{document}
\title{\textbf{Domain Walls in the Heisenberg-Ising Spin-$\fr{1}{2}$ Chain}}
\author{Axel Saenz\footnote{Department of Mathematics, Oregon State University, Corvallis, OR 97331,
email: saenzroa@oregonstate.edu } \and Craig A.~Tracy\footnote{Department of Mathematics,
UC Davis, Davis, CA 95616, email: tracy@math.ucdavis.edu} \and Harold Widom\footnote{Department of Mathematics,
UC Santa Cruz, Santa Cruz, CA 95064}}
\date{}
\maketitle

\begin{abstract}
In this paper we obtain formulas for the distribution of the left-most up-spin in the Heisenberg-Ising spin-1/2 chain with anisotropy parameter $\Delta$, also known as the XXZ spin-1/2 chain, on the one-dimensional lattice $\mathbb{Z}$ with domain wall initial conditions. We use the Bethe Ansatz to solve the Schr\"odinger equation and a recent antisymmetrization identity of Cantini, Colomo, and Pronko \cite{CCP} to simplify the marginal distribution of the left-most up-spin. In the $\Delta=0$ case, the distribution $F_2$ arises. In the $\Delta \neq 0$ case, we propose a conjectural series expansion type formula based on a saddle point analysis. The conjectural formula turns out to be a Fredholm series expansion in the $\Delta \rightarrow 0$ limit and recovers the result for $\Delta = 0$.
\end{abstract}

\section{Introduction}
We consider the dynamics of the Heisenberg-Ising spin-1/2 chain   with  anisotropy parameter $\Delta$, also known as the XXZ spin-1/2 chain,
on the one-dimensional lattice $\bZ$ with domain wall initial conditions.  We start with an initial state of $N$ up-spins at the
sites $\{1,2,\ldots, N\}$ in a sea of down-spins; and by utilizing ideas from coordinate Bethe Ansatz \cite{Bethe, Gaudin, Suth, YY} to solve the
Schr\"odinger equation, we find the quantum state $\Psi_N(t)$ at time $t$ is
\[ \Psi_N(t)=\sum_X \psi_N(X,t) e_X, \]
where the sum is over all $X=\{x_1<x_2<\cdots<x_N\}$ and $e_X$ denotes the state with up-spins at $X$. Alternatively we can view a spin up at site $x_j$ as a particle and a spin down as an empty lattice site or hole. The ``Bethe-coordinates'' $\psi_N(X,t)$ are given below in Theorem 1.\footnote{Since the Hamiltonian $H_{\Delta}$ 
of the Heisenberg-Ising model is a 
(non-unitary) similarity transformation of the Markov generator of the ASEP \cite{Gwa-Spohn}, the results in \cite{TW1} give immediately
the Bethe coordinates of Theorem 1 once
an identification of parameters is made (see Section 3.3).} They have the
interpretation that $\left\vert\psi_N(X;t)\right\vert^2$ is the probability the system is in state $X$ at time $t$. Observe that the $\psi_N(X,t)$ have the standard Bethe Ansatz structure
as a sum over the permutation group $\cS_N$; where now, each term in the summand is an $N$-dimensional contour integral.
\par
\subsection{One-Point Functions}
If $X_1(t)$ denotes the position of the left-most particle at time $t$, then
\[ \bP_N(X_1(t) = x)=\sum_{X,x_1=x}\left\vert\psi_N(X,t)\right\vert^2 \]
where the sum is over all $X=\{x_1=x<x_2<\cdots<x_N\}$.  In ASEP the analogous quantity involves a \textit{single sum} over $\cS_N$ where as now
we have a \textit{double sum} over $\cS_N$.
In \cite{TW1} an identity involving the sum over the permutation group\footnote{See equation (1.6) in \cite{TW1}.} was used to reduce the
sum to a single  $N$-dimensional integral.  Cantini, Colomo, and Pronko \cite{CCP} have  generalized the single sum permutation
identity to a double sum permutation identity, which also generalize to the (spin) Hall-Littlewood functions \cite{P, WZJ}. Employing this new identity reduces the expression for $\bP_N(X_1(t)=x)$ to
a single $2N$-dimensional integral whose integrand involves the famous Izergin-Korepin determinant \cite{Iz,Ko}.  The resulting expression
is given in Theorem 2. This part of the paper overlaps the recent work of J.~M.~St\'ephan \cite{St1, St2}.
\par For the special case $\Delta=0$,  the analysis simplifies considerably.  Using Toeplitz operators and their determinants,  we show the $N\rightarrow\infty$ limit
can be taken resulting in the representation
\[ \lim_{N\ra\infty}\bP_N(X_1(t)\ge x )=\det(I-L)\]
where $L$ is an integral operator whose kernel is the  \textit{discrete Bessel kernel} \cite{Bo,BOO, Jo}. (See also Chapter 8 in
\cite{BDS}.) This makes connections to the distribution
of the length of the longest
increasing subsequence in a random permutation \cite{BDJ, BDS}.  See Theorem 3 below.  From this identification it follows that
\[ \lim_{t\ra\iy}\lim_{N\ra\iy}\bP_N\left(\fr{X_1(t)+2t}{t^{1/3}}\ge -y\right)=F_2(y)\]
where $F_2$ is the $\textrm{TW}_2$ distribution \cite{TW0a, TW0b}. This last result appears to be well-known in the physics literature since
the case $\Delta=0$ is
reducible to a ``free fermion'' model \cite{ER, SZ, St1, VSDH}.

\subsection{Contour Deformations and a Conjecture}

Taking the contour integral functions for the one-point function to the infinite time statistics is another major challenge. In the case of the ASEP, this was achieved by Tracy-Widom \cite{TW1a} by deforming the contours to obtain a Fredholm determinant. Then, in a later work by the same authors \cite{TW2}, the Fredholm determinant was further analyzed by deforming the kernels to obtain the Tracy-Widom distribution. We also deform the contour integrals for our one-point function, in Section \ref{s:contour_deformations}, to obtain a type of series expansion.

\textsc{Theorem 4.} Let $X_1(t)$ be the location of the left-most particle in the Heisenber-Ising spin-1/2 chain with $N$ particles, initial conditions $Y = (y_1< y_2 < \cdots < y_N )$, and $\D \in \mathbb{R}$ so that $\D \neq 0$. Then, $\mathbb{P}(X_1(t) \geq x)$ is equal to
\be
\sum_{n=0}^N \sum_{\tau \in \mathcal{T}_n}\oint_{\cC_R}\cdots \oint_{\cC_R} \oint_{\cC_{R'}} \cdots \oint_{\cC_{R'}} I_N(\x, \z; \tau) f(\x, \z; \tau)\, \left(\prod_{j \in J} d\z_j \right)\, d^N \x\hspace{5ex}
\ee
where the integrand is given by (\ref{e:integrands}), the summation is take over the set of maps $\mathcal{T}_n$ given by (\ref{e:maps_cont}), and the contours $\cC_R$ and $\cC_{R'}$ are circles centered at zero with radii $R, R' > 0$ that satisfy the following inequalities $\max\{2|\D|^{-1} , 2(1+2|\D|)\} < R < \max\{4|\D|^{-1}, 4(1 + 2|\D|)\} < R'/2$.

We expect this series expansion to to give rise to a series expansion of a Fredholm determinant in the infinite time limit. In fact, we may deform the contours in the previous formula to the steepest descent to contours in an effort to obtain the infinite time limit by a saddle point analysis. The result is given by our \textsc{Conjecture 8.1}. Aside from technical details of certain bounds and approximations, there are some terms that we still can't control after the saddle point analysis. Recent results \cite{BK,CdLV, St1, St1a}, based on numerical, hydrodynamic and analytical arguments are inconclusive in the appropriate scaling, i.e.~$t^{1/2}$ versus $t^{1/3}$, for the fluctuations of the one-point function in the infinite time limit. Based on our conjecture, we expect the location of the left-most particle to be at $-2t$ with fluctuations on the order of $t^{1/3}$ but the limiting distribution is still unclear.

\section{XXZ Quantum Spin-$\hf$ Hamiltonian}\label{s:XXZ}
The definition of the quantum spin chain Hamiltonian on the infinite lattice $\bZ$ requires some  explanation since there is the problem of making sense of  infinite tensor products in the construction of a Hilbert space of states. The general construction uses the Gelfand-Naimark-Segal (GNS) construction; but in the case considered here, there is an elementary treatment 
\cite{NSS} which we now describe.
\par
 Let $\cH_{0}=\bC$.  For each positive integer $N$ we define
\[ \cX_N:=\left\{X=\{x_1,\ld,x_N\}\in\bZ^N: x_1<\cdots<x_N\right\} \]
and 
\[ \cH_N:=\ell^2(\cX_N).\]
The Hilbert space of states is 
\[ \cH:=\bigoplus_{N=0}^\iy \cH_N.\]
The normalized state $\Om=1\in\cH_0$ is the ground state of all spins down.  In physicists' notation
\[ \Om= \vert\cd\da\cd\da\cd\da\cd\rangle.\]
Given $N\in\bZ^{+}$ and $X=\{x_1,\ld,x_N\}\in\cX_N$, define $e_X\in\cH_N$ by
\[ e_X(Y)=\dl_{X,Y}.\]
The set $\{e_X\}_{X\in\cX_N}$ defines a natural orthonormal basis of $\cH_N$. The physical interpretation of
$e_X$ is the state with up spins at $x_1<\cd<x_N$ in a sea of down spins:
\[ e_X=\vert \cd \underset{x_1}{\ua}\cd\underset{x_2}{\ua}\cd\underset{x_N}{\ua}\cd\rangle. \]
This is a model
of a quantum lattice gas (see, for example, \S6.1.6 of \cite{Suth}).  We will frequently use this particle interpretation.\par
We  introduce the \textit{Pauli operators} $\s_j^\a$, $j\in\bZ$, $\a=3,\pm$.
\bae
 \s_j^3 e_X&=&\left\{\begin{array}{ll}\hspace{1.5ex} \> e_X&\textrm{if}\>\>\>j\in X=\{x_1,\ld,x_N\},\\
								-\> e_X & \textrm{otherwise}.\end{array}\right.\\ 
\s_j^{+} e_X&=&\left\{\begin{array}{ll} 0 &\textrm{if}\>\>\> j\in \{x_1,\ld,x_N\},\\
e_{X^{+}} & \textrm{where}\>\>\> X^+=\{x_1,\ld,x_k, j, x_{k+1},\ld, x_N\}, \> x_k<j<x_{k+1}\end{array}\right.\\
 \s_j^{-} e_X&=&\left\{\begin{array}{ll} 0 &\textrm{if}\>\>\> j\notin X=\{x_1,\ld,x_N\}\\
e_{X^{-}}& \textrm{where}\>\>\>X^{-}=\{x_1,\ld,x_{k-1},x_{k+1},\ld,x_N\},  \>\>\> j=x_k\end{array}\right.\eae

In words, $\s_j^{+}:\cH_N\ra\cH_{N+1}$  acts as the identity except at the site $j$ where
it takes $\da \leadsto \ua$ and annihilates a $\ua$ state.  Similarly, $\s^{-}_j:\cH_{N}\ra\cH_{N-1}$ acts as the identity except at the site $j$
where it takes $\ua\leadsto \da$ and annihilates a $\da$ state.  By definition 
$\s_j^3\Om=-\Om$, $\s_j^{-}\Om=0$ and $\s_j^{+}\Om=e_{\{j\}}$.    We also recall the 
Pauli  operators  $\s_j^1=\s_j^{+}+\s_j^{-}$ and $\s_j^2=-i\s_j^{+}+i\s_j^{-}$.
Define
 \[ h_{j,j+1}=\fr{1}{2}\left(\s_j^1\s_{j+1}^1
+\s_{j}^2\s_{j+1}^2+\D(\s_j^3\s_j^3-1)\right) =\s_j^{+}\s_{j+1}^{-}+\s_j^{-}\s_{j+1}^{+}+\fr{\D}{2}(\s_j^3\s_{j+1}^3-1)\]
and
\be H_{XXZ}=\sum_{j\in\bZ} h_{j,j+1} \label{XXZ}.\ee
The operator $H_{XXZ}$ is the \textit{Heisenberg-Ising spin-$\fr{1}{2}$ chain Hamiltonian}; or more briefly, the  $XXZ$ spin Hamiltonian.
It's clear from the above definitions that $H_{XXZ}:\cH_N\ra\cH_N$.  Since the number of particles is conserved under the dynamics of $H_{XXZ}$, we can work in a sector $\cH_N$.
\par
A state $\Psi_N=\Psi_N(t)\in\cH_N$ can be represented by
\be  \Psi_N(t)=\sum_{X\in\cX_N} \psi_N(X,t) e_X.\label{BigPsi}\ee
The initial condition is $\Psi(0)=e_Y$, $Y=\{y_1,\ld, y_N\}\in\cX_N$,  so that $\psi_N(X;0)=\dl_{X,Y}$.  The dynamics is determined by the Schr\"odinger equation
\be \mi\, \fr{\pl\Psi_N}{\pl t}=H_{XXZ}\Psi_N.\label{SchroEqn}\ee
The Hamiltonian $H_{XXZ}$ is self-adjoint and so by Stone's theorem there
exists a unitary operator \newline $U=\exp(-\mi t H_{XXZ})$ such that $\Psi_N(t)=U(t)\Psi_N(0)$.  We have
\[ \langle \Psi_N(t),\Psi_N(t)\rangle =\sum_{X\in \cX_N} \left\vert \psi_N(X;t)\right\vert^2=1.\]

The goal is to describe the dynamics $\Psi_{DW}(t)$ starting from the \textit{domain wall} (DW) initial state 
\[ e_{\bN}=\vert\cd \da\da\underset{0}{\da}\underset{1}{\ua}\ua\ua\cd\rangle.\]
One immediately sees the difficulty in that $e_{\bN}$ is not an element of $\cH_N$ for any 
$N$.\footnote{Presumably,
one could construct a domain wall Hilbert space $\cH_{DW}$ by replacing the state $\Omega$ by
$e_{\bN}$.  Unfortunately, we do not know how to proceed with a Bethe Ansatz solution in this space.}  If
$X_m(t)$ denotes the position of the $m$th particle on the left, we define
\[ \bP_{\bN}(X_m(t)=x)=\lim_{N\ra\iy} \bP_{\{1,\ldots,N\}}(X_m(t)=x).\]

\section{Bethe Ansatz Solution $\Psi_N(t)$}\label{s:bethe_ansatz}
This section closely follows \cite{TW1, YY}.  We first note that 
\bae
h_{j,j+1}\lvert \cd \underset{j}{\ua}\underset{j+1}{\ua}\cd\rangle&=&0,\label{h1}\\
h_{j,j+1}\lvert \cd \underset{j}{\da}\underset{j+1}{\da}\cd\rangle&=&0,\label{h2}\\
h_{j,j+1}\lvert \cd \underset{j}{\ua}\underset{j+1}{\da}\cd\rangle&=&
-\D \lvert \cd \underset{j}{\ua}\underset{j+1}{\da}\cd\rangle +
\lvert \cd\underset{j}{\da}\underset{j+1}{\ua}\cd\rangle,\label{h3}\\
 h_{j,j+1}\lvert \cd \underset{j}{\da}\underset{j+1}{\ua}\cd\rangle&=&
-\D \lvert \cd \underset{j}{\da}\underset{j+1}{\ua}\cd\rangle +
\lvert \cd\underset{j}{\ua}\underset{j+1}{\da}\cd\rangle.\label{h4}
\eae

\subsection{$N=1$}
Let $\Psi_1(t)=\sum_{x_1}\psi_1(x_1;t) e_{\{x_1\}}$, then
\ba H_{XXZ}\Psi_1(t)&=&\sum_{j} h_{j,j+1}\sum_{x_1}\psi_1(x_1;t) e_{\{x_1\}}\\
&=&  \sum_j\psi_1(j;t) h_{j,j+1} e_{\{j\}}+\psi_1(j+1;t)h_{j,j+1}e_{\{j+1\}}\\
&=& \sum_j \psi_1(j;t)\left[ -\D e_{\{j\}}+e_{\{j+1\}}\right]+\psi_1(j+1;t)
\left[-\D e_{\{j+1\}}+e_{\{j\}}\right]\\
&=&\sum_j \left[ -2\D\psi_1(j;t)+\psi_1(j-1;t)+\psi_1(j+1;t)\right] e_{\{j\}}
\ea
Thus the coordinates $\psi_1(x;t)$ must satisfy
\[ i\fr{\pl\psi_1(x;t)}{\pl t}= \psi_1(x-1;t)+\psi_1(x+1;t)-2\D \psi_1(x;t),\>\>x\in\bZ, t\ge 0,\]
with initial condition
\[ \psi_1(x;0)=\dl_{x,y}.\]
The solution is
\be \psi_1(x;t)=\fr{1}{2\pi i}\,\int_{\cC_r}\x^{x-y-1} e^{-i t\ve(\x)}d\x\, =e^{2i\D t} (-i)^{x-y} \, J_{x-y}(2t)\label{psi1}\ee
where
\[ \ve(\x)=\x+\fr{1}{\x}-2\D\]
and $J_\nu(z)$ is the Bessel function of order $\nu$.
\par
It is obvious probabilistically, and is easily verified analytically, that
\[ \sum_{x\in\bZ} \vert \psi_1(x;t)\vert^2 =1 \]
for all $y\in\bZ$.  We also have (setting $y=0$)
\[ \sum_{x\in \bZ} x \vert\psi_1(x;t)\vert^2=0,\>\> \sum_{x\in\bZ} x^2\vert\psi_1(x;t)\vert^2=t,
\>\>\textrm{and}\>\>\sum_{x\in\bZ} x^4\vert\psi_1(x;t)\vert^2=t^2+3 t^4.
\]
\subsection{$N=2$}
For $N=2$ we set
\[ \Psi_2(t)=\sum_{x_1<x_2}\psi_2(x_1,x_2;t) e_{\{x_1,x_2\}}\]
For $x_2>x_1+1$ the action of $H_{XXZ}$ on $\Psi_2$ is the same as above in each coordinate $x_i$:
\be i \fr{\pl \psi_2(x_1,x_2;t)}{\pl t}=\psi_2(x_1-1,x_2;t)+\psi_2(x_1+1;x_2;t)+\psi_2(x_1,x_2-1;t)+\psi_2(x_1,x_2+1)-4\D\psi_2(x_1,x_2;t).
\label{n2a}\ee
If $x_2=x_1+1$ then due to (\ref{h1}) and (\ref{h2}) there are terms missing with the result that
\be i \fr{\pl \psi_2(x_1,x_2;t)}{\pl t}=\psi_2(x_1-1,x_2;t)+\psi_2(x_1,x_2+1)-2\D\psi_2(x_1,x_2;t).\label{n2b}\ee
We now require (\ref{n2a}) to hold for all $(x_1,x_2)\in\bZ^2$ and to satisfy the boundary condition
\[ \psi_2(x_1,x_1;t)+\psi_2(x_1+1,x_1+1;t)-2\D\psi_2(x_1,x_1+1;t)=0,\>\> x_1\in\bZ.\]
If this last boundary condition is satisfied then in region $\cX_2$ equation (\ref{n2b}) is satisfied.

Define \cite{YY} (the Yang-Yang $S$-matrix)
\be S_{21}(\x_2,\x_1)=-\fr{1+\x_1\x_2-2\D\x_2}{1+\x_1\x_2-2\D\x_1},\quad  \x_1,\x_2\in\bC.\label{YY}\ee
With this choice of $S$,
\[ \left\{ \x_1^{x_1-y_1-1}\x_2^{x_2-y_2-1}+S_{21}(\x_2,\x_1) \x_2^{x_1-y_2-1}
\x_1^{x_2-y_1-1}\right\} e^{-it(\ve(\x_1)+\ve(\x_2))}\]
satisfies (\ref{n2a}) and the boundary condition (\ref{n2b}); however it does not satisfy the initial condition.
We take\footnote{ Here and later all differentials $d\x$ and $d\z$ incorporate the factor $(2\pi \mi)^{-1}$.} ``linear combinations''
\be \psi_2(x_1,x_2;t)=\int_{\cC_r}\int_{\cC_r}\left\{ \x_1^{x_1-y_1-1}\x_2^{x_2-y_2-1}+S_{21}(\x_2,\x_1) \x_2^{x_1-y_2-1}
\x_1^{x_2-y_1-1}\right\} e^{-it(\ve(\x_1)+\ve(\x_2))}\, d\x_1 d\x_2.\label{psi2}\ee
A residue calculation shows that in the physical region\footnote{$x_1<x_2$ and $y_1<y_2$} the initial condition $\psi(x_1,x_2;0)=\delta_{x_1,y_1}\delta_{x_2,y_2}$ will be  satisfied if the radius of the circle $\cC_r$,
centered at the origin, is sufficiently small so that the only singularities of the integrand
lying inside $\cC_r$ are those at 0.
\subsection{General $N$}

The generator of the finite $N$ asymmetric simple exclusion process (ASEP) is a similarity transformation 
(\textit{not} a unitary transformation!) of the Heisenberg-Ising Hamiltonian.  Because of this the Schr\"odinger equation (\ref{SchroEqn})
for the quantum spin chain is essentially identical to the master equation (Kolmogorov forward equation) for the Markov process ASEP assuming
the identification of parameters
\[ \x_i= \x_i^\prime /\sqrt{\t} ,\>\>\t=\fr{p}{q},\>\>2\D=\fr{1}{\sq{pq}},\>\>
S^{XXZ}_{\b\a}(\x_{\b},\x_{\a})=S_{\b\a}^{ASEP}(\x_{\b}',\x_{\a}'),\>\>
\ve^{XXZ}(\x)=\fr{1}{\sqrt{pq}}\ve^{ASEP}(\x').\]
\par
Thus given the  ASEP result \cite{TW1, TW3} and the above identifications, we have
\par
\textsc{Theorem 1.}  
For $\s\in\cS_N$, define
\be A_\s(\x)=\prod\left\{S_{\b\a}(\x_\b,\x_\a): \{\b,\a\}\>\> \textrm{is an inversion in}\>\> \s\right\}, \label{e:A_coeff}\ee
then the solution to (\ref{SchroEqn}) satisfying the initial condition $\psi_N(X;0)=\delta_{X,Y}$ is
\be \psi_N(X;t)=\sum_{\s\in\cS_N}\int_{\cC_r}\cd\int_{\cC_r} A_\s(\x) \prod_i\x_{\s(i)}^{x_i} \prod_i \left(\xi_i^{-y_{i}-1}\,
\me^{-\mi t \ve(\xi_i)}\right)
 d\x_1\cd d\x_N\label{psi_N}\ee
where $\cC_r$ is a circle centered at zero with radius $r$ so small that all the poles of $A_\s$ lie outside of
$\cC_r$.

Additionally, we have a contour integral formula with large contours instead of small contours as above in \textsc{Theorem 1}. Below, we will use a combination of the small and large contour formulas.

\textsc{Theorem 1a.}  
For $\s\in\cS_N$, define
\[ A_\s(\x)=\prod\left\{S_{\b\a}(\x_\b,\x_\a): \{\b,\a\}\>\> \textrm{is an inversion in}\>\> \s\right\},\]
then the solution to (\ref{SchroEqn}) satisfying the initial condition $\psi_N(X;0)=\delta_{X,Y}$ is
\be \psi_N(X;t)=\sum_{\s\in\cS_N}\int_{\cC_R}\cd\int_{\cC_R} A_\s(\x) \prod_i\x_{\s(i)}^{x_i} \prod_i \left(\xi_i^{-y_{i}-1}\,
\me^{-\mi t \ve(\xi_i)}\right)
 d\x_1\cd d\x_N\label{psi_N_large}\ee
where $\cC_R$ is a circle centered at zero with radius $R$ so large that all the poles of $A_\s$ lie inside of $\cC_R$.

The proof of this statement is an adaptation of the arguments in \cite{TW1} that give the proof of \textsc{Theorem 1}. For completeness, we give the proof of \textsc{Theorem 1a} in \Cref{s:appendix_a}.

 \section{Probability $\cP_Y(x,m;t)$}\label{s:one_point}
If the initial state is $e_Y\in\cH_N$, $Y\in\cX_N$, then
at time $t$ the system is in state $\Psi_N(t)=\sum_{X\in\cX_N} \psi_N(X;t) e_X$ where
$\psi_N(X;t)$ is given by (\ref{psi_N}) or (\ref{psi_N_large}).  The
quantity
\[ \left\vert\langle e_X,\Psi_N(t)\rangle\right\vert^2=\left\vert\psi_N  (X;t)\right\vert^2 ,\>\> X\in\cX_N,\]
is the probability that the system is in state $e_X$ at time $t$. 
\par
Denote by $\cP_Y(x,m;t)$ the probability that at time $t$ the state has the $m$th particle from the left at position $x$
given initially the state is $Y$. 
Let $X=\{x_1,x_2,\ldots, x_N\}\in\cX_N$, $1\le m\le N$, and define the projection operator
\be P_{x,m}e_X=\left\{\begin{array}{cc}e_X&\textrm{if}\>\> x_m=x,\\
0& \textrm{otherwise}.\end{array}\right.\ee 
Then the outcome of the measurement yielding  ``the $m$th spin from the left is at position $x$ at
time $t$'' is that the system is now in state
\[\Psi_N(x,m;t):= P_{x,m}\Psi_N(t)=\sum_{\substack{X\in\cX_N \\x_m=x}}  \psi_N(X;t) e_X.\]
Thus the probability of this outcome is
\be \cP_Y(x,m;t):=\langle\Psi_N(x,m;t),\Psi_N(x,m;t)\rangle=
\sum_{\substack{X\in\cX_N \\x_m=x}}\left\vert\psi_N(X;t)\right\vert^2. \label{1ptProb}\ee

\subsection{Distribution of left-most particle}
We now restrict to the case $m=1$, i.e.\ $\cP_Y(x,1;t)$.
Let
\[ x_1=x,\> x_2=x+v_1,\>\ldots, x_N=x+v_1+v_2+\cdots+v_{N-1},\> v_i\ge 1,\]
and  note that  $\overline{\Psi_N(x;t)}=\Psi_N(x;-t)$. Then, using (\ref{psi_N}) for $\Psi(x;t)$ and (\ref{psi_N_large}) for $\Psi(x;-t)$ with $R \, r < 1$, followed by performing the
geometric sums (since $R\,r <1$, the summations may be brought inside) 
\ba
\cP_Y(x,1;t)&=&\sum_{\substack{X\in\cX_N \\x_1=x}} \psi_N(X;t)\psi_N(X;-t)\\
&=&\sum_{\s,\mu\in \cS_N}\int_{\cC_R}\cdots\int_{\cC_r} \sum_{v_i\ge 1} A_\s(\x) A_\mu(\z)\,
(\x_{\s(2)}\z_{\mu(2)})^{v_1} (\x_{\s(3)}\z_{\mu(3)})^{v_1+v_2}\cdots 
(\x_{\s(N)}\z_{\mu(N)})^{v_1+\cdots+v_{N-1}}\\
&&\times
\prod_j (\x_j\z_j)^{x-y_j-1} \me^{-\mi t(\ve(\x_j)-\ve(\z_j))}\,d\z_1\cdots d\z_N d\x_1\cdots d\x_N\\
&\hspace{-15ex} =&\hspace{-10ex}\sum_{\s,\mu\in \cS_N}\int_{\cC_R}\cdots\int_{\cC_r}  A_\s(\x) A_\mu(\z)
\fr{\x_{\s(2)}\z_{\mu(2)} \x_{\s(3)}^2\z_{\mu(3)}^2\cdots
 \x_{\s(N)}^{N-1} \z_{\mu(N)}^{N-1}}{(1-\x_{\s(2)}\z_{\mu(2)}\cdots \x_{\s(N)}\z_{\mu(N)})
(1-\x_{\s(3)}\z_{\mu(3)}\cdots \x_{\s(N)}\z_{\mu(N)})\cdots (1-\x_{\s(N)}\z_{\mu(N)})}\\
&&\times\prod_j (\x_j\z_j)^{x-y_j-1} \me^{-\mi t(\ve(\x_j)-\ve(\z_j))}\,d\z_1\cdots d\z_N d\x_1\cdots d\x_N 
\ea
In the formulas above, we have $2N$ contour integrals with the contour $\cC_r$ for the first $N$ contours and the contours $\cC_R$ for the following $N$ contours.
Now, at the analogous step in ASEP, an identity\footnote{See (1.6) in \cite{TW1}.}  was
derived that simplified the sum over $\cS_N$
 resulting in a single multidimensional integral.\footnote{See Theorem 3.1 in \cite{TW1}.}
 Now we have a \textit{double sum} over $\cS_N$ and we need a new identity.  Fortunately
 such an identity has been discovered by Cantini, Colomo, and Pronko \cite{CCP}.
 Let
 \be d(x,y):=\fr{1}{(1-x \,y)(x+y -2\Delta\,  x\, y)}\>\>\textrm{and}\>\> 
D_N(\x,\z)=\det\left(d(\x_i,\z_j)\vert_{1\le i,j\le N}\right),\label{IK}\ee
then
\bae
\sum_{\s,\mu\in\cS_N} A_\s(\x) A_{\mu}(\z) \fr{\x_{\s(2)}\z_{\mu(2)} \x_{\s(3)}^2\z_{\mu(3)}^2\cdots
 \x_{\s(N)}^{N-1} \z_{\mu(N)}^{N-1}}{(1-\x_{\s(2)}\z_{\mu(2)}\cdots \x_{\s(N)}\z_{\mu(N)})
(1-\x_{\s(3)}\z_{\mu(3)}\cdots \x_{\s(N)}\z_{\mu(N)})\cdots (1-\x_{\s(N)}\z_{\mu(N)})}\nonumber\\
=\fr{(1-\prod_j \x_j\z_j) \, 
\prod_{i,j=1}^N(\x_i+\z_j-2\Delta \x_i\z_j)}{
\prod_{i<j} (1+\x_i\x_j-2 \Delta \x_i) (1+\z_i\z_j-2 \Delta \z_i)}\>D_N(\x,\z)&&\label{1ptIden}
\eae
Remarks:\vspace{-4ex}
\begin{itemize}
\item  The identity (\ref{1ptIden}) is Proposition 6 of \cite{CCP} (with a change of notation). The identity (\ref{1ptIden}) also appears in a more general setting of (spin) Hall-Littlewood functions in \cite{P, WZJ}, which specializes to the ASEP case as shown in Corollary 7.1 in \cite{P}.
\item
  In Appendix B of \cite{CCP}, the authors show that
 (\ref{1ptIden}) reduces to (1.6) of \cite{TW1} 
in the limit $\x_j\ra \sqrt{\fr{q}{p}}\,\x_j$ and $\z_j\ra \sqrt{\fr{p}{q}}$.
\item The determinant $D_N(\x,\z)$ ``is nothing but the well-known \textit{Izergin-Korepin
determinant} \cite{Iz, Ko} in disguise'' \cite{W}. 
\end{itemize}
\par
We thus have
\vspace{-2ex}
\be \hspace{-5ex}\cP_Y(x,1;t)=\int_{\cC_R}\cdots\int_{\cC_r}\fr{(1-\prod_j \x_j\z_j) \, 
\prod_{i,j=1}^N(\x_i+\z_j-2\Delta \x_i\z_j)}{
\prod_{i<j} (1+\x_i\x_j-2 \Delta \x_i) (1+\z_i\z_j-2 \Delta \z_i)}\>D_N(\x,\z)
\prod_j (\x_j\z_j)^{x-y_j-1} \me^{-\mi t(\ve(\x_j)-\ve(\z_j))}\,d^N\z d^N\x
\label{detRep}\ee
The factor $(1-\prod_j \x_j\z_j)$ is eliminated if we consider 
\be \cF_N(x,t):=\bP_Y(X_1(t)\ge x)=\sum_{n=x}^\iy\cP_Y(n,1;t)\label{cFn}\ee
\par

From \cite{CCP}
\be  \prod_{1\le j,k\le N} (\x_j+\z_k-2 \D \x_j \z_k)\,\cdot\, D_N(\x,\z)
=\fr{\D_N(\x)\D_N(\z)}{\prod_{j,k} (1-\x_j\z_k)} Q_{N}(\x,\z) \label{Wfn} \ee
where $Q_N$ is a ``polynomial of degree $N-1$ in each variable, separately symmetric under permutations
of the variables within each set'' \cite{CCP, W}.\footnote{For example
\ba Q_1(\x,\z)&=&1,\\
Q_2(\x,\z)&=&4 \Delta ^2 \zeta _1 \zeta _2 \xi _1
   \xi _2-2 \Delta  \zeta _1 \zeta
   _2 \xi _1-2 \Delta  \zeta _1
   \zeta _2 \xi _2-2 \Delta  \zeta
   _1 \xi _1 \xi _2-2 \Delta  \zeta
   _2 \xi _1 \xi _2+\zeta _1 \zeta
   _2 \xi _1 \xi _2+\zeta _1 \zeta
   _2+\xi _1 \xi _2+1,
   \ea
   $Q_3$ in expanded form has 459 terms, and $Q_4$ has 60,820 terms.} Here $\D_N(\x)$ is the Vandermonde product $\prod_{1\le j<k\le N}(\x_k-\x_j)$ (not to
be confused with the constant $\D$).
   It's useful to define
\[ U(\x,\x'):=\fr{1+\x\x'-2\D\x}{\x'-\x} .\]

   The identity (\ref{1ptIden}) can be rewritten as
   \bae \sum_{\s,\mu}\prod_{i<j} U(\x_{\s(i)},\x_{\s(j)}) U(\z_{\mu(i)},\z_{\mu(i)})\,
   \fr{\x_{\s(2)}\z_{\mu(2)} \x_{\s(3)}^2\z_{\mu(3)}^2\cdots
 \x_{\s(N)}^{N-1} \z_{\mu(N)}^{N-1}}{(1-\x_{\s(2)}\z_{\mu(2)}\cdots \x_{\s(N)}\z_{\mu(N)})
(1-\x_{\s(3)}\z_{\mu(3)}\cdots \x_{\s(N)}\z_{\mu(N)})\cdots (1-\x_{\s(N)}\z_{\mu(N)})}\nonumber \\
&\hspace{-160ex}=&\hspace{-80ex} \fr{1-\prod_j\x_j \z_j}{\prod_{j,k}(1-\x_j\z_k)}\, Q_N(\x,\z)
\label{idenU}
\eae
The close relationship of (\ref{idenU}) to (1.6) of \cite{TW1} (see also Identity $1_L$ in \cite{TW4}) is now clearer.
We have proved
\par
\textsc{Theorem 2.}
$\cF_N(x,t)=\bP_Y(X_1(t)\ge x)$ equals
\be
\int_{\cC_R}\cdots\int_{\cC_r} \fr{\prod_{j,k}(\x_j+\z_k-2\D\x_j\z_k)}{\prod_{j<k}(1+\x_j\x_k-2\D\x_j)(1+\z_j\z_k-2\D\z_j)}\,
D_N(\x,\z)
\prod_j (\x_j \z_j)^{x-y_j-1} \me^{-\mi t(\ve(\x_j)-\ve(\z_j))}\, d^N\z \,d^N\x\hspace{5ex}\label{Fprob1}\ee
\be = \int_{\cC_R}\cdots\int_{\cC_r} \fr{\D_N(\x)\D_N(\z)}{\prod_{j<k}(1+\x_j\x_k-2\D\x_j)
(1+\z_j\z_k-2\D\z_j)}\,
\fr{Q_N(\x,\z)}{\prod_{j,k}(1-\x_j\z_k)}
\prod_j (\x_j \z_j)^{x-y_j-1} \me^{-\mi t(\ve(\x_j)-\ve(\z_j))}\, d^N\z \,d^N\x\hspace{5ex}\label{FprobQ}
\ee
where $\cC_r$ (resp.~$\cC_R$) is a circle centered at zero with radius $r$ (resp.~$R$) so small (resp.~large) that all the poles of the integrand except for the the poles at the origin (resp.~infinity) lie outside $\cC_r$ (resp.~inside $\cC_R$) and $R\,r <1$.
\section{Special case $\D=0$}\label{s:delta_zero}
When $\D=0$, (\ref{FprobQ}) reduces to
\bae
\cF_N(x,t)\big\vert_{\D=0}&=&
    \int_{\cC_R}\cdots\int_{\cC_r} \fr{\D_N(\x)\D_N(\z)}{\prod_{j,k}(1-\x_j\z_k)}
\prod_j (\x_j \z_j)^{x-y_j-1} \me^{-\mi t(\ve(\x_j)-\ve(\z_1))}\, d^N\z \,d^N\x\hspace{5ex}
\label{FprobQ2}\\
&=& 
\int_{\cC_R}\cdots\int_{\cC_r} \det\left(\fr{1}{1-\x_j\z_k}\right)
\prod_j (\x_j \z_j)^{x-y_j-1} \me^{-\mi t(\ve(\x_j)-\ve(\z_1))}\, d^N\z \,d^N\x\hspace{5ex}
\label{FprobQ3}
\eae
since
\[ \lim_{\D\ra 0}\fr{Q_N(\x,\z)}{\prod_{j<k}(1+\x_j\x_k-2\D\x_j)(1+\z_j\z_k-2\D\z_j)}=1.\]
and
\[ \det\left(\fr{1}{1-\x_j\z_k}\right)=\fr{\D_N(\x)\D_N(\z)}{\prod_{j,k} (1-\x_j\z_k)}.\]
 More directly since $A_\s\big\vert_{\Delta=0}=\sgn(\s)$,  we use the identity
    \[
    \sum_{\s,\mu}\sgn(\s)\sgn(\mu) \fr{\x_{\s(2)}\z_{\mu(2)} \x_{\s(3)}^2\z_{\mu(3)}^2\cdots
 \x_{\s(N)}^{N-1} \z_{\mu(N)}^{N-1}}{(1-\x_{\s(2)}\z_{\mu(2)}\cdots \x_{\s(N)}\z_{\mu(N)})
(1-\x_{\s(3)}\z_{\mu(3)}\cdots \x_{\s(N)}\z_{\mu(N)})\cdots (1-\x_{\s(N)}\z_{\mu(N)})}\]
   \be =\det\left(\fr{1}{1-\x_j\z_k}\right)
    \ee
    \par
    \subsection{Fredholm determinant representation}
Define
\[ \phi_j(\x)=\x^{x-y_j-1}\me^{-\mi t\ve(\x)},\>\> 
\psi_j(\z)=\z^{x-y_j-1}\me^{\mi t\ve(\z)}\]
and
\be K(j,k)=\fr{\phi_j(\x_j) \psi_k(\z_k)}{1-\x_j\z_k}\label{K}\ee
Thus
\be 
\cF_N(x,t)\big\vert_{\D=0}=
    \int_{\cC_R}\cdots\int_{\cC_r} \det(K)\, d^N\z \,d^N\x= \int_{\cC_r}\cdots\int_{\cC_r} \det(K)\, d^N\z \,d^N\x\hspace{5ex}\label{Kdet}\ee
    \par
    For the second identity, we deformed the contours from $\cC_R$ to $\cC_r$ for all the $\zeta$-variables. When we deform the contours, we don't cross any poles since the poles, given by $1- \xi_j \zeta_k=0$, are located outside of the contour $\cC_R$ since we have taken $R\, r <1$. Additionally, note that the variable $\x_j$ appears only in row $j$ and $\z_k$ appears only in column $k$. It follows that the multiple integral
    is gotten by integrating each $K(j,k)$ with respect to $\x_j$, $\z_k$.
    Therefore the multiple integral (\ref{Kdet}) equals the determinant with $j,k$ entry
    \[ K_N(j,k)=\int_{\cC_r}\int_{\cC_r} \fr{\phi_j(\x)\psi_k(\z)}{1-\x\z}\,d\z d\x \]
    \par
    We consider step initial condition, so that $y_j=j$. In preparation
    for taking the limit as $N\ra\iy$, we make the replacements
    $j\ra j+1$, $k\ra k+1$, so that the indicies run for $0$ to $N-1$  rather than
    $1$ to $N$. Then, in preparation for eventual steepest descent, we make
    the substitutions $\x\ra\mi\,\x$, $\z\ra\z/\mi$.  Aside from the factor $\me^{\mi\pi(j-k)/2}$, which will not
    affect the determinant, the kernel becomes
        \[ L_N(j,k)=\int_{\cC_r}\int_{\cC_r}\fr{\x^{x-j-2}\,\z^{x-k-2}}{1-\x\z} 
        \me^{t (\theta(\x)+\theta(\z))}\, d\z d\x,\]
    where we have set $\theta(\x)=\x-1/\x$.  We write the above as
    \[ \sum_{\ell=0}^\iy\int_{\cC_r}\int_{\cC_r} \x^{x-j+\ell-2}\z^{x-k+\ell-2}\,\me^{t(\theta(\x)+\theta(\z))}\,d\z d\x.\]
    \par
    We may take all integrations over the unit circle $\cC_1$ and in the $\z$-integral make the substitution
    $\z\ra1/\z$.  We obtain
    \[ L_N(j,k)=\sum_{\ell=0}^\iy \int_{\cC_1}\int_{\cC_1} \x^{x-j+\ell-2}\z^{-x+k-\ell}\,\me^{t(\theta(\x)-\theta(\z))}\, d\z d\x.\]
    In Toeplitz terms this is the operator
    \[ P_N T(a) T(a^{-1}) P_N,\]
    where $P_N$ is the projection from $\ell^2(\bZ^+)$\footnote{$\bZ^+$ denotes
    the set of nonnegative integers.} to $\ell^2([0,\ldots,N-1])$ and where $a$ is the symbol
    \[ a(\x)=\x^{x-1} \me^{t\,\theta(\x)}.\]
    It it known (see, e.g.\ \S5.1 in \cite{BS})  that $T(a) T(a^{-1})$ is of the form $I$+trace class and so $\det(K_N)$ has
    the limit $\det(T(a) T(a^{-1}))$ on $\ell^2(\bZ^+)$.\footnote{One can show that for $x>1$ the determinant
    of the product is zero.}
    \par
    By a well-known identity, $T(a) T(a^{-1})=I-H(a) H(\tilde{a}^{-1})$, where $H(a)$ denotes
    the Hankel operator and $\tilde{a}(\x)=a(\x^{-1})$.  In this case $\tilde{a}=a^{-1}$ and the square
    of $H(a)$ has kernel\footnote{Recall that the $i,j$-entry of $H(f)$ is $f_{i+j+1}=\int \x^{-i-j-2} f(\x)\,d\x$.}
    \[ L(j,k)=\sum_{\ell=0}^\iy \int_{\cC_1}\int_{\cC_1} \x^{x-j-\ell-3} \z^{x-k-\ell-3} \,\me^{t(\theta(\x)+\theta(\z))}
    \, d\z d\x,\]
    and we are interested in $\det(I-L)$.   The substitutions $\x\ra1/\x$, $\z\ra 1/\z$ give
    \be L(j,k)=\sum_{\ell=0}^\iy \int_{\cC_1}\int_{\cC_1}\x^{-x+j+\ell+1} \z^{-x+k+\ell+1}\,\me^{-t(\theta(\x)+\theta(\z))}\, d\z d\x.\label{Lexpand}\ee
    If we take our integrals over $\cC_r$ and sum we obtain
    \be L(j,k)=\int_{\cC_r}\int_{\cC_r}\fr{\x^{-x+j+1}\z^{-x+k+1}\,\me^{-t(\theta(\x)+\theta(\z))}}{1-\x\z}\ d\z d\x
    \label{L}\ee
    The kernel $L(j,k)$ is known as the \textit{discrete Bessel kernel} \cite{Bo} (see also Chapter 8 in \cite{BDS}) due to the following representation.
    Using the Bessel generating function
    \[\exp(t\theta(\x))=\sum_{n=-\iy}^\iy \x^n J_n(2t) \] 
    in (\ref{Lexpand}) and the identity, $\nu\neq \mu$,
    \be \sum_{n=0}^\iy J_{\nu+n}(t) J_{\mu+n}(t)= \fr{t}{2(\nu-\mu)}\left[J_{\nu-1}(t) J_{\mu}(t)-J_{\nu}(t) 
    J_{\mu-1}( t)\right]\label{besselIden}\ee
we find
\[ L(j,k) =t\, \fr{J_{j-x+1}(2t) J_{k-x+2}(2t)-J_{j-x+2}(2t) J_{k-x+1}(2 t)}{j-k}\]
For $j=k$ one lets $\mu\ra\nu$ in (\ref{besselIden}) to find
\[L(j,j)=\sum_{n=0}^\iy J_{\nu+n}(2t)^2 =t\,\left[ J_\nu(2t) \fr{\pl J_\mu}{\pl\mu}\big\vert_{\mu=\nu-1} -J_{\nu-1}(2t) 
\fr{\pl J_\mu}{\pl \mu}\big\vert_{\mu=\nu-1}\right],\>\nu=-x+j+1.\] 
For $x\le 1$ and domain wall initial condition $Y=\bN$,  we have the Toeplitz representation 
\be \bP_{\bN}(X_1(t)\ge x)\big\vert_{\D=0}=\det(I-L)_{\ell^2(\{1-x, 2-x, \dots\})}=\me^{-t^2} \det\left(I_{j-k}(2 t)\right)
 \Big\vert_{j,k=0,\ldots,-x}
\label{BOiden}\ee
where the last equality\footnote{$I_\nu(z)$ is the modified Bessel function of order $\nu$.} was proved in \cite{BOk}.
  \par
  If $\cL(t)$ denotes the length of the longest increasing subsequence of a random permutation of size $\cN$ where
  $\cN$ is a Poisson random variable with parameter $t^2$, then \cite{BDJ, BDS, G}
  \[ \bP(\cL(t)\le n) =\me^{-t^2} \det(I_{j-k}(2t))_{j,k=0,\ldots,n-1}\]
  
\textsc{Theorem 3.} For $x \leq 1$ and domain wall initial conditions $Y = \mathbb{N}$, we have
\be
\bP_{\bN}(X_1(t)\ge x)\big\vert_{\D=0}=\bP(\cL(t)\le 1-x)
\ee
where $\cL(t)$ denotes the length of the longest increasing subsequence of a random permutation of size $\cN$ so that $\cN$ is a Poisson random variable with parameter $t^2$.

\subsection{Asymptotics}
From the classic work of Baik, Deift, and Johnasson \cite{BDJ} (see also Chapter 9 in \cite{BDS}), we know that the limiting distribution of
$\cL(t)$ is
\be \lim_{t\ra\iy}\bP\left(\fr{\cL(t)-2 t}{t^{1/3}}\le x\right) =F_2(x)\ee
where $F_2$ is the $\beta=2$  TW distribution  \cite{TW0a, TW0b}.
In the present problem, $\D=0$,  we can therefore conclude 
that the left-most particle for domain wall initial condition $Y=\bN$  has the limiting distribution
\be \lim_{t\ra\iy} \bP\left(\fr{X_1(t)+2t}{t^{1/3}}\ge -y\right)= F_2(y).\label{TWF2}\ee

\section{Steepest descent curve}\label{s:steepest_descent}

\subsection{Spectral functions}
We introduce a pair of functions 
\begin{equation}\label{e:spectral_function}
        G(\xi) = x \log \xi - i t (\xi + \xi^{-1}), \quad H(\zeta) = -x \log \zeta - i t (\zeta + \zeta^{-1}),
\end{equation}
which we call the \emph{spectral functions}. Note that the spectral functions appear in the integrand of the formula for $\mathcal{F}_N(x,t)$ given by (\ref{Fprob2}). In particular, we have
\begin{equation}
    (\xi_j \zeta_j)^{x} e^{- i t (\varepsilon(\xi_j) - \varepsilon(\zeta_j))} = \exp \left \{ G(\xi_j) - H(\zeta_j)\right\}.
\end{equation}
In the following, we will deform the contours in the contour integral formula for $\mathcal{F}_N$ given by (\ref{Fprob2}) so that the real part of the difference of the spectral function is negative, $\mathrm{Re}(G- H) <0$. Thus, making $\mathcal{F}_N$ suitable for asymptotic analysis.

\subsection{Critical points}
The steepest descent contours in the contour integral formula $\mathcal{F}_N$ given by (\ref{Fprob1}) are determined by the critical points of the spectral functions. We have
\begin{equation}
        G'(\xi)  = \frac{- i t \xi^2 + x \xi + i t}{\xi^2}, \quad H'(\zeta) = \frac{-i t \zeta^2 - x \zeta + it }{\zeta^2}.
\end{equation}
so that the critical points are given by
\begin{equation}
    \xi = \frac{x \pm \sqrt{x^2 -4 t^2}}{2 i t},\quad \zeta = \frac{-x \pm \sqrt{x^2 -4 t^2}}{2 i t}.
\end{equation}
Note that each function, $G$ and $H$, has a double critical point when $x = \pm 2 t$ and the critical point are
\begin{equation}
    \xi_0  = \begin{cases} - i , \quad x= 2 t \\ i , \quad x = -2t \end{cases}, \quad \zeta_0  = \begin{cases}  i , \quad x= 2 t \\ -i , \quad x = -2t \end{cases},
\end{equation}
respectively. Physically, we expect the point $x=-2t$ to correspond to the left-edge of the up-spins and the point $x= 2t$ to correspond to the right-edge of the up-spins. Thus, we restrict our attention to the critical point given by $x= -2t$ and take $(\xi_0, \zeta_0)= (i, -i)$.

\subsection{Steepest descent curve: local}

\begin{figure}[h!]
    \centering
    \includegraphics[scale=0.3]{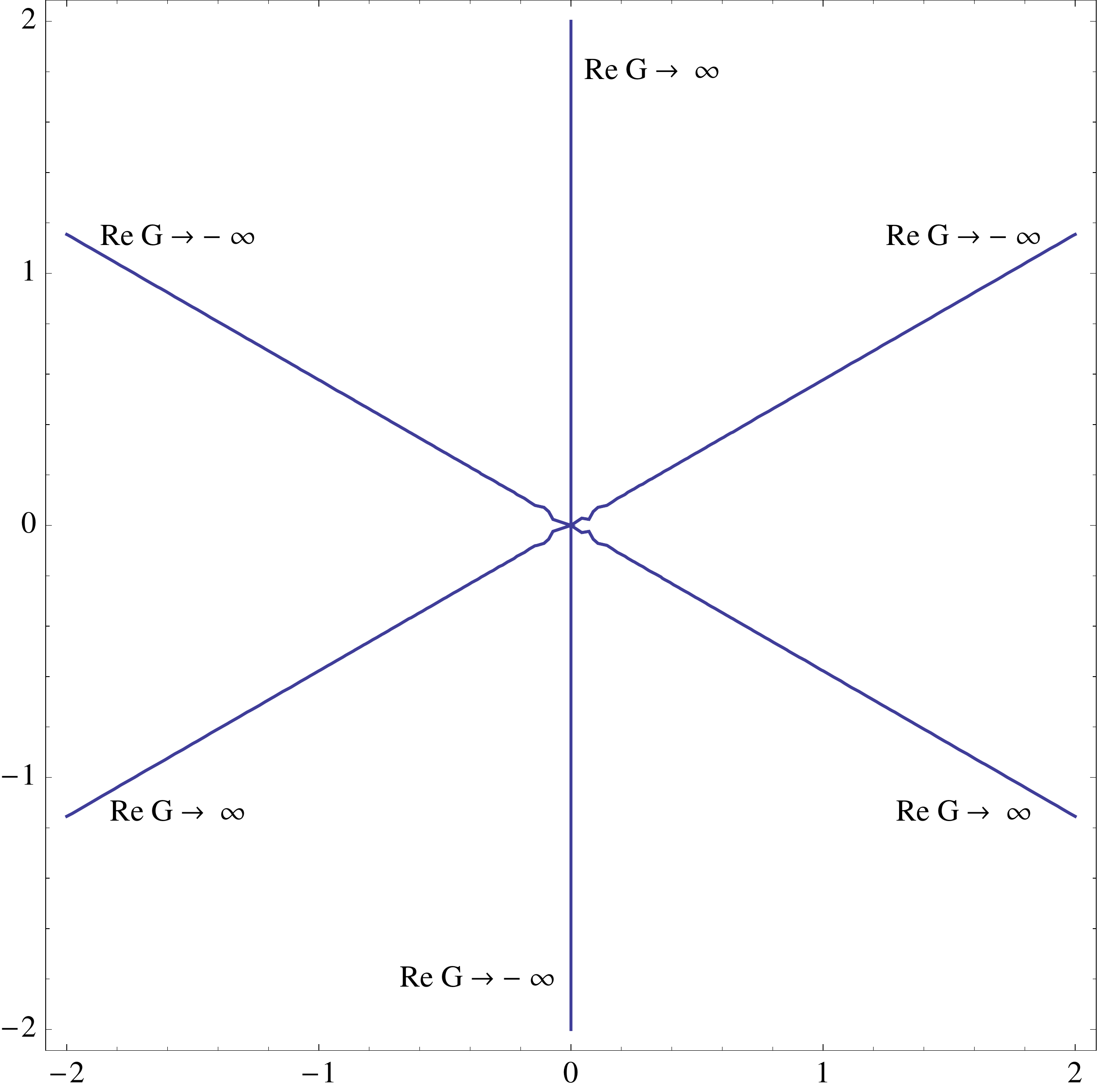}
    \caption{The curves show the directions of large increase of magnitude for the real part of $G$ (i.e.~steepest descent) near the critical point, which is centered to be at the origin. We label the curves so that the sign of $Re\, G$ is clear.}
    \label{f:local_SD}
\end{figure}

The steepest descent curve depends locally on the value of the third derivative for the double critical points. In particular, we have
\begin{equation}
    G^{(3)}(- 2t) = 2 i t, \quad H^{(3)}(- 2 t) = 2 i t.
\end{equation}
This means that the steepest descent curves for both spectral functions are locally the same. That is,
\begin{equation}\label{e:Taylor}
    G(\xi) - G(\xi_0) = \frac{i t}{3} (\xi - \xi_0)^3 + \mathcal{O}((\xi-\xi_0)^4), \quad H(\zeta) - H(\zeta_0) = \frac{i t}{3} (\zeta - \zeta_0)^3 + \mathcal{O}((\zeta-\zeta_0)^4) 
\end{equation}
A diagram for the local steepest descent curve is give in Figure \ref{f:local_SD}.

\subsection{Steepest descent curve: global}

\begin{figure}
     \centering
     \begin{subfigure}[h!]{0.45\textwidth}
         \centering
         \includegraphics[width=\textwidth]{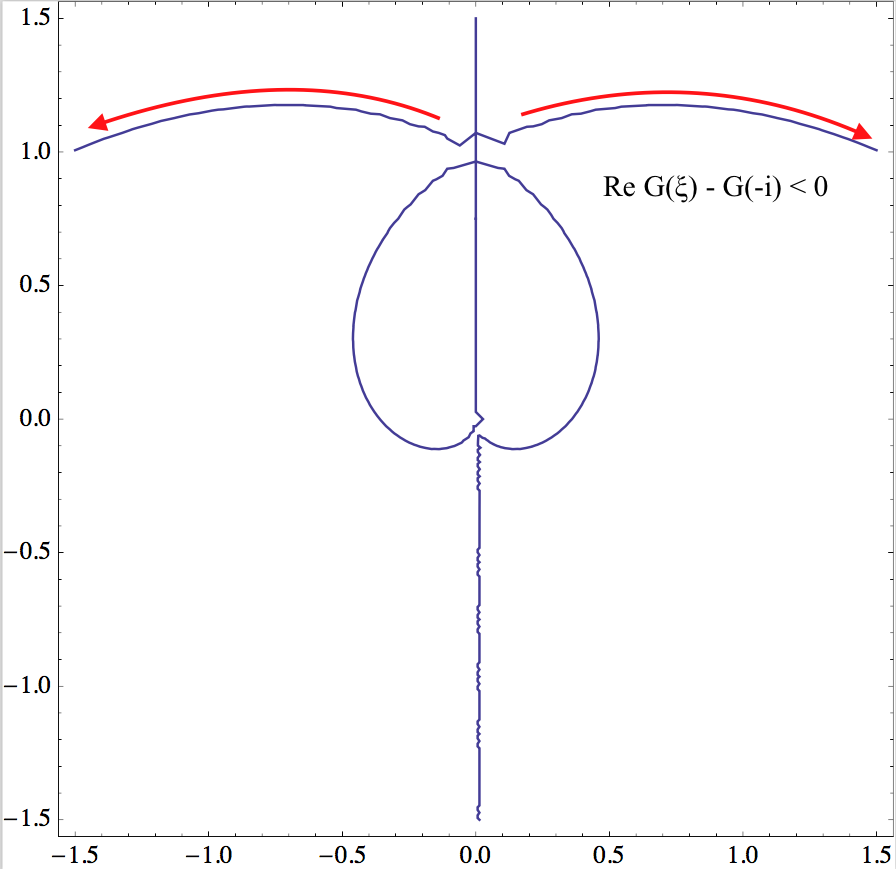}
         \caption{}
     \end{subfigure}
     \hfill
     \begin{subfigure}[h!]{0.45\textwidth}
         \centering
        \includegraphics[width=\textwidth]{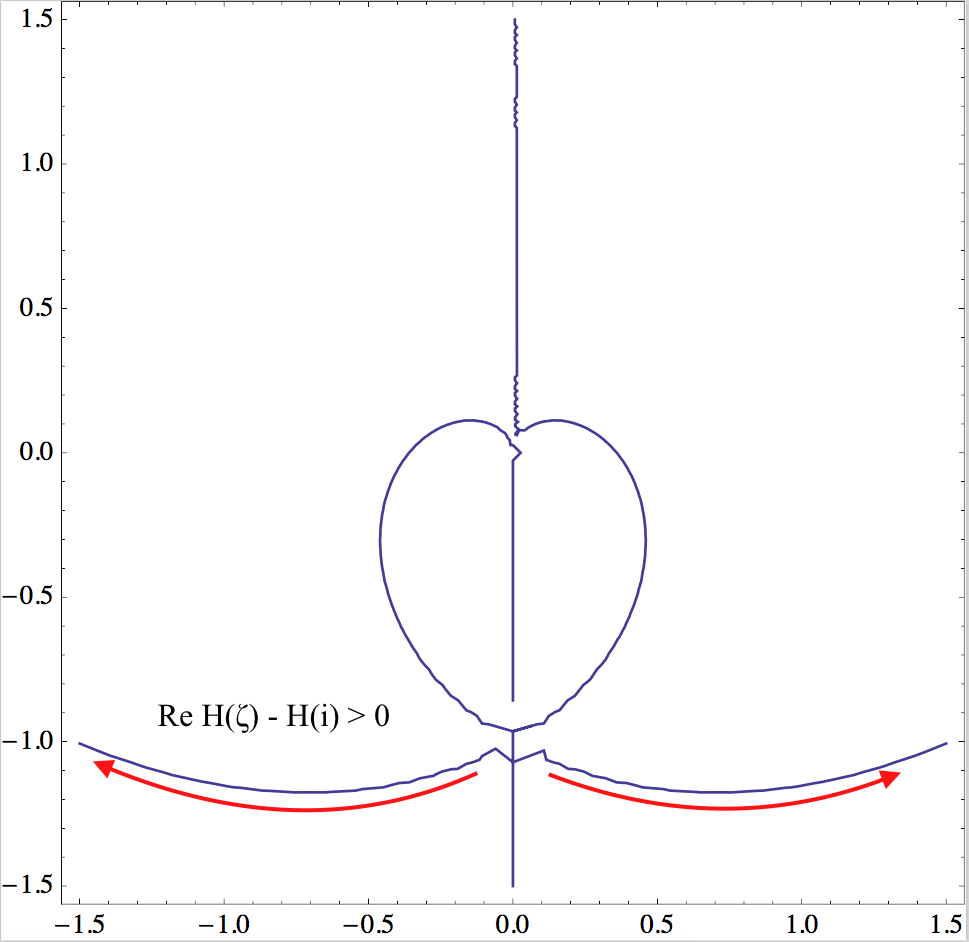}
         \caption{}
     \end{subfigure}
        \caption{The level curves for the spectral functions, $G$ and $H$, with $\xi_0 = i$, $\zeta_0=-i$, and $x = -2 t$. The component corresponding to the steepest descent path is denoted by red arrows with the direction along the steepest descent.}
        \label{f:global_SD_2}
\end{figure}

Globally, we obtain the steepest descent curve by an implicit equation through the Cauchy-Riemann equations for holomorphic functions. In short, we have that the direction of greatest change in the real part of a holomorphic function is the same direction of no change in the imaginary part of a the same holomorphic function (away from singular points). Therefore, the steepest descent curve is given by the level curves of the imaginary part of the spectral functions:
\begin{equation}
    \{\xi \in \mathbb{C} : Im\, G(\xi) = Im\, G(\xi_0)\}, \quad \{\zeta \in \mathbb{C} : Im\, H(\zeta) = Im\, H(\zeta_0)\}
\end{equation}
with $(\xi_0 , \zeta_0) = ( i, -i)$. The level curves will have multiple components, but we are only interested in the components that (i) are closed curves in $\mathbb{C} \cup \{ \infty\}$ and (ii) $Re \, G(\xi) - G(\xi_0) \leq 0$ and $Re \, H(\xi) - H(\xi_0) \geq 0$.\par

The steepest descent curves are components of the curves given by
\begin{equation}\label{e:steepest_descent}
\begin{split}
    &2 \tan^{-1}\left(\frac{y}{x} \right) + \frac{x^3 + x y^2 +x}{x^2 +y^2} = \begin{cases}  \pi, &\quad x>0 \\ -\pi, &\quad x<0 \end{cases},\\
    & 2 \tan^{-1}\left(\frac{v}{u} \right) - \frac{u^3 + u v^2 +u}{u^2 +v^2} =  \begin{cases}  -\pi, &\quad u>0 \\ \pi, &\quad u<0 \end{cases}
\end{split}
\end{equation}
with $\xi = x + i y$ and $\zeta = u + i v$. We plot these curves using \emph{Mathematica}, see Figure \ref{f:global_SD_2}. The components of the level curves corresponding to the steepest descent path is determined by the local behaviour determined in Figure \ref{f:local_SD}.

\subsection{Steep descent curve}

We want a friendly version of the steepest descent curves given implicitly by (\ref{e:steepest_descent}) or, rather, a more explicit version. We introduce the \emph{steep descent contours} given by three segments on three regions: in the region near the critical points, we take straight lines coming emanating from the critical point at angles $\pm\pi/6$ and $\pm5\pi/6$; in an intermediate region, we take horizontal lines emanating from the end points of the straight lines in region near the critical point; in the region far away from the critical point, we take a segment of the circle $\cC_R$ that connects with the horizontal lines. We use these contours so that we may explicitly determine the location of the poles when we deform to these \emph{steep descent contours}. Although these contours don't follow the path of steepest descent for the real part of the spectral function, we show below that we still have the main property that $Re\,\{ G(\xi) - G(\xi_0)\} \leq 0$ and $Re \, \{H(\zeta) - H(\zeta_0)\} \geq0$ along these steep descent contours.

\begin{figure}
     \centering
     \begin{subfigure}[h]{0.45\textwidth}
         \centering
         \includegraphics[width=\textwidth]{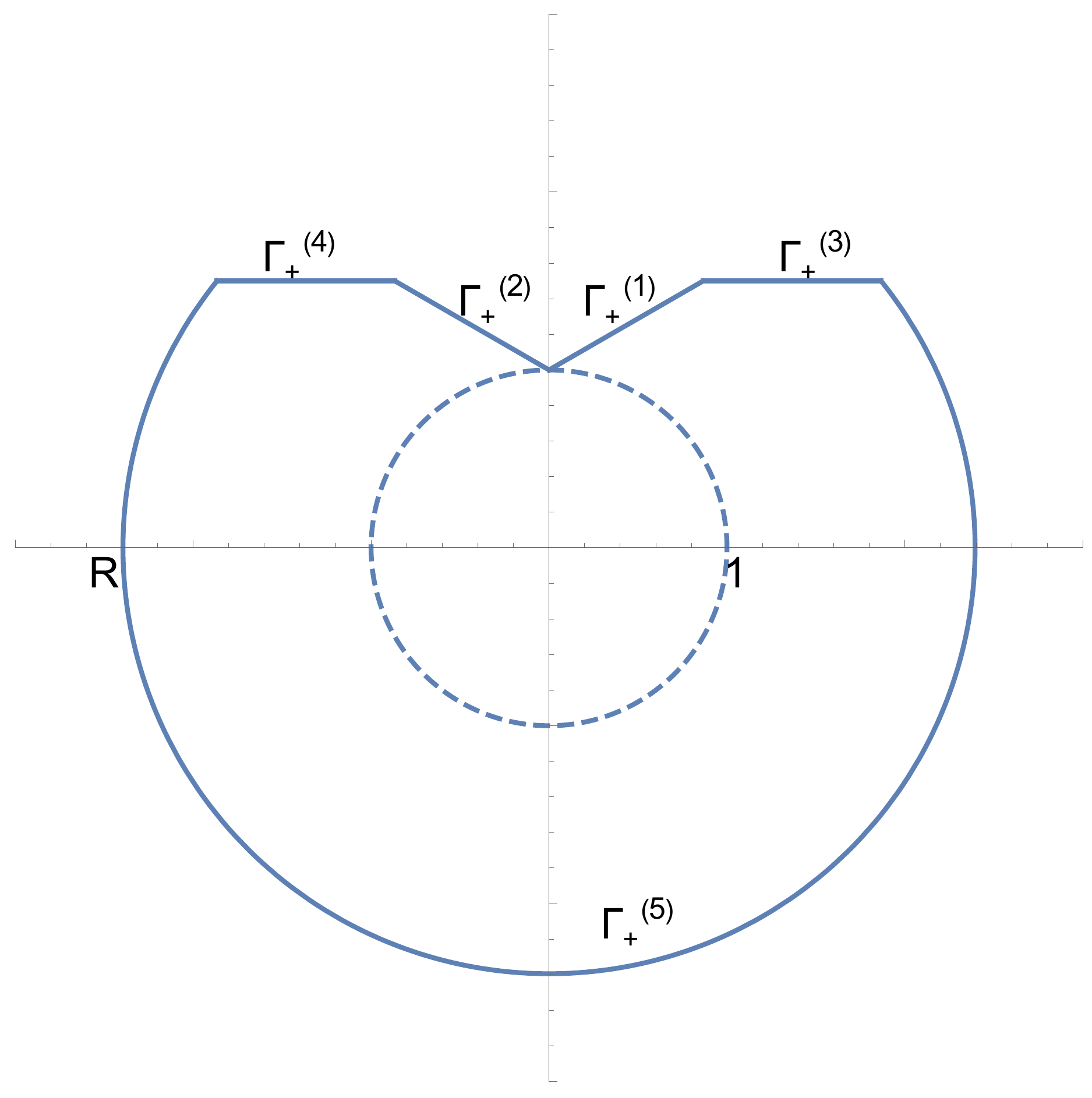}
         \caption{}
     \end{subfigure}
     \hfill
     \begin{subfigure}[h]{0.45\textwidth}
         \centering
        \includegraphics[width=\textwidth]{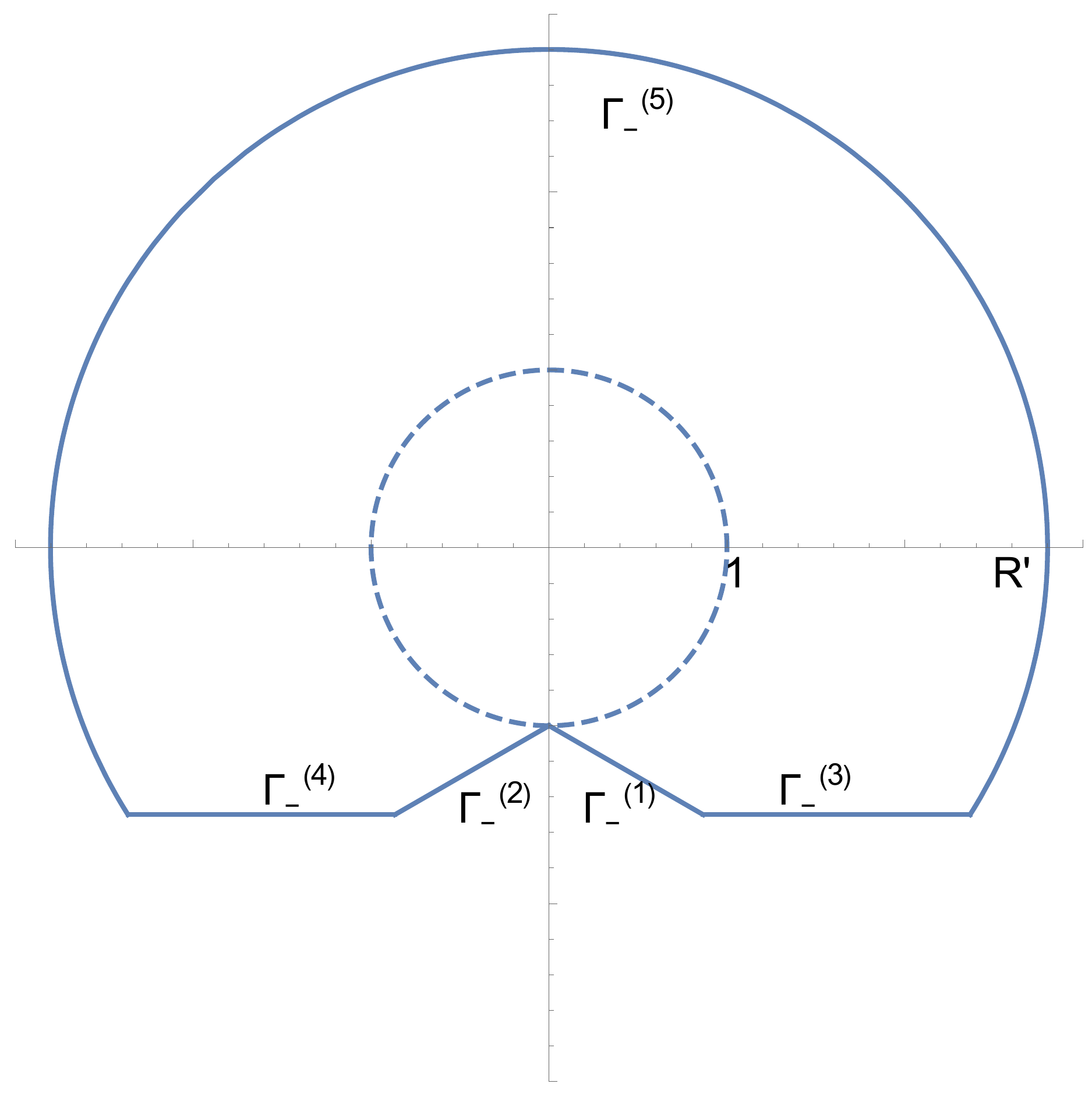}
         \caption{}
     \end{subfigure}
        \caption{The components of the $\Gamma_{\pm}$ contours.}
        \label{f:global_SD_3}
\end{figure}

We now give a precise definition for the steep descent contours. We give a piece-wise description based on the proximity to the critical points. Let $\cB(z, r)$ be a ball centered at $z\in \mathbb{C}$ of radius $r>0$ and $\cB(z,r)^c$ be its complement. Then, we take the components 
\begin{equation}\label{e:contour_parts}
    \begin{split}
    \Gamma_{\pm}^{(1)}&=\{\pm i + x  e^{ \pm\pi i/6} \mid 0 \leq x \} \cap \cB(\pm i, 1),\\
    \Gamma_{\pm}^{(2)}&=\{\pm i + x e^{ \pm 5\pi i/6} \mid 0 \leq x \} \cap \cB(\pm i, 1)\\
    \Gamma_{\pm}^{(3)}&=\{\pm i +  e^{\pm \pi i/6} +x \mid 0 \leq x \} \cap \cB(\pm i, 1)^c \cap \cB(0, R_{\pm}),\\
    \Gamma_{\pm}^{(4)}&=\{\pm i + 1 e^{\pm 5 \pi i/6} -x \mid 0 \leq x \} \cap \cB(\pm i, 1)^c \cap \cB(0, R_{\pm})\\
    \Gamma_{\pm}^{(5)}&= \cC_R \cap \{z\in \mathbb{C} \mid Im \, \{z\} \leq (\pm1)Im\, \{\pm i + e ^{\pm \pi i /6}\} \}
    \end{split},
\end{equation}
with radii $R_{\pm} > \sqrt{3}$. The bound on the radii is chosen so that the horizontal segments of the contours are non-trivial. Then, the steep descent contours are given by
\begin{equation}\label{e:steep_contours}
    \Gamma_{k} = \Gamma_k^{(1)} \cup \Gamma_k^{(2)} \cup \Gamma_k^{(3)} \cup \Gamma_k^{(4)} \cup \Gamma_k^{(5)}
\end{equation}
for $k=\pm$. See Figure \ref{f:global_SD_3}.

\textsc{Lemma 6.1} Let $x=-2t$ and take the contours $\Gamma_k$, $k=\pm$, given by (\ref{e:contour_parts}) and (\ref{e:steep_contours}). Additionally, take $t^{-\alpha} \leq T \ll 1$ with $1/4 < \alpha < 1/3$. Then, we have
\begin{equation}\label{e:bound1}
    Re\,\{ G(\xi) - G(\xi_0)\} \leq 0, \quad Re \, \{H(\zeta) - H(\zeta_0)\} \geq 0 
\end{equation}
if $\xi \in \Gamma_{+}$ and $\zeta \in \Gamma_{-}$. Moreover, if $\xi \in \Gamma_{+} \cap \cB(\mi, t^{-\alpha})^c$ and $\zeta \in \Gamma_{-} \cap \cB(-\mi, t^{-\alpha})^c$, we have
\begin{equation}\label{e:bound2}
    Re\,\{ G(\xi) - G(\xi_0)\} < - c_1(T)\, t^{1- 3\alpha}, \quad Re \, \{H(\zeta) - H(\zeta_0)\} > c_2(T)\, t^{1-3\alpha}, 
\end{equation}
for some constants $c_1(T), c_2(T)>0$ that depend only on $T$.

\begin{proof}
We prove the bounds by showing that derivative of the real part of the functions are monotone along the different segments of the contours $\Gamma_{\pm}$ as parameterized in (\ref{e:contour_parts}). Since $G(\xi) - G(\xi_0) = 0$ for $\xi = \xi_0$ and $H(\zeta) - H(\zeta_0) = 0$ for $\zeta = \zeta_0$, the first bounds (\ref{e:bound1}) then follow by monotonicity. Moreover, since the real part of the functions are monotone, we establish the bounds (\ref{e:bound2}) by bounding the real part of the functions on the boundary of the segment $\Gamma_{\pm} \cap \cB(\pm \mi, t^{-\alpha})$. 

The arguments for both functions are the same, except for some negative signs here and there. So, we focus solely on the case for the $G$ function. Additionally, the arguments are fairly routine and standard. So, we just sketch the main idea needed for the bounds.

Take $\xi \in \Gamma_{+}^{(1)} \cup \Gamma_{+}^{(2)}$. In this case, we have $\xi =\mi + x e^{\pi i /6}$ or $\xi  = \mi +x  e^{5\pi i /6}$, with $0 \leq x \leq 1$ since $\Gamma_{+}^{(1)} \cup \Gamma_{+}^{(2)} \subset \cB(i, 1)$. Then, we may write the real part of the $G$ function explicitly and show that it is monotone by taking its derivative. For instance, we have
\begin{equation}
    \frac{d}{dx} Re\,\{G(i +  x e^{\pi i /6} ) - G(i)\} = \frac{t}{2} \left(1 - \fr{2 +4 x}{1 +x + x^2} + \fr{1 + 4x +x^2}{(1 + x+x^2)^2} \right).
\end{equation}
One may now check that the derivative is zero when $x=0$ and negative if $0 < x < 1 + \sqrt{3}$. Thus, the bound (\ref{e:bound1}) follows for this segment.

Take $\xi \in \Gamma_{+}^{(3)} \cup \Gamma_{+}^{(4)}$. In this case, we have $\xi = i +  e^{\pi i /6} + x$ or $\xi  = i+ e^{5\pi i /6} -x$, with $x$ non-negative and bounded since $\Gamma_{+}^{(3)} \cup \Gamma_{+}^{(4)} \subset \cB(0, R_{+})$. Then, we may write the real part of the $G$ function explicitly and show that it is monotone by taking its derivative. For instance, we have
\begin{equation}
    \frac{d}{dx} Re\,\{G(i +   e^{\pi i /6} + x) - G(i)\} = -t \left(1 - \frac{3(\sqrt{3}  + 2 x)}{2(3 + \sqrt{3}\, x + x^2)^2} \right).
\end{equation}
Form this, one may show that the derivative is strictly negative for all $x \geq 0$. The bound (\ref{e:bound1}) follows for this segment.

Take $\xi \in \Gamma_{+}^{(5)}$. In this case, we have $\xi = R_{+}\, e^{i \theta}$, with $- \pi/2 \leq \theta \leq \phi_1 < \pi/2$ and $\pi/2 < \phi_2 \leq \theta \leq 3\pi/2$ for some constants $\phi_1$ and $\phi_2$ since $\Gamma_{+}^{(5)} \subset \{z \in \mathbb{C} \mid Im\, \{z\}  \leq Im\, \{ i +  e^{\pi i /6}\} \}$. In this case, we have
\begin{equation}
    Re\,\{G(\xi) - G(i)\} = -2 t\,\log R_{+} + t(R_{+}+R_{+}^{-1}) \sin \theta.
\end{equation}
Since $R_{+}>1$, one may then show that this function is monotone on $\theta$ for each of the segments $- \pi/2 \leq \theta \leq \phi_1 < \pi/2$ and $\pi/2 < \ph_2 \leq \theta \leq 3\pi/2$. The bound (\ref{e:bound1}) follows for this segment.

The bound (\ref{e:bound2}), now that we have established that the function is monotone along all the segments of the contours, follows by evaluating the function on the boundary of the segment $\Gamma_{+} \cap \cB(i, t^{-\alpha})$. That is, we evaluate the function at the points $\x = \x_0 + t^{-\alpha} \, e^{ \pi \mi/6}$ and $\x = \x_0 + t^{-\alpha} e^{5 \pi \mi/6 }$. In particular, we use the Taylor expansion
\begin{equation}
    G(\xi) - G(\xi_0) = -\frac{1}{3}x^3 t^{1- 3 \alpha} + \mathcal{O}(t^{1-4\alpha})
\end{equation}
to approximate the function at the desired points. Since $t^{-\alpha}< T \ll 1$, we obtain the bound (\ref{e:bound2}). 
\end{proof}

\section{Contour Deformations}\label{s:contour_deformations}

\subsection{Small to Large Contour deformations}

We deform the contours in the probability function for the left-most particle given by (\ref{Fprob1}). In particular, we deform the contours $\cC_r$, for the $\zeta$-variables, to some contour $\cC_{R'}$ with a large radius $R'>0$. Let
\begin{equation}\label{e:deformation_contour}
    \Omega(\xi) :=  \mathcal{C}^{(0)} \cup -\cC^{(1)} \cup -\cC^{(2)} \cup \cdots \cup -\cC^{(N)}  
\end{equation}
be the union of $(N+1)$ circles so that $-\cC^{(j)}$, for $j=1, \dots, N$, is a negatively oriented circle centered at $\xi_j^{-1}$ with radius $ r' >0$ and $\cC^{(0)}$ is a positively oriented circle centered at the origin with radius $R'>0$. We give precise conditions on the radii in the statement of \textsc{Lemma 7.1} below. Then, as we deform the $\cC_r$ contour, we will encounter poles at $\zeta_i = \xi_j^{-1}$ for $i, j =1, \dots, N$. As a result, we obtain the contour $\Omega(\xi)$ when we deform the contour $\cC_r$ to $\cC_{R'}$. This result and the proof for the contour deformations, given by \textsc{Lemma 7.1} below, is similar to the contour deformation in \cite{CDP}.

\textsc{Lemma 7.1}
For $\Delta \neq 0$, $\cF_N(x,t)=\bP_Y(X_1(t)\ge x)$ equals
\be
\int_{\cC_R}\cdots\int_{\Omega(\xi)} \fr{\prod_{j,k}(\x_j+\z_k-2\D\x_j\z_k)}{\prod_{j<k}(1+\x_j\x_k-2\D\x_j)(1+\z_j\z_k-2\D\z_j)}\,
D_N(\x,\z)
\prod_j (\x_j \z_j)^{x-y_j-1} \me^{-\mi t(\ve(\x_j)-\ve(\z_j))}\, d^N\z \,d^N\x\hspace{5ex}\label{Fprob2}\ee
where the contour $\cC_R$ for the $\xi$-variables is a circle centered at zero with radius $R>0$ and the contour $\Omega(\xi)$ for the $\zeta$-variables is given by (\ref{e:deformation_contour}) with radii $R' >0$ and $r' = 1/(2 R)$, so that the radii satisfy the following inequalities $\max\{2|\D|^{-1}, 2(1+2|\D|)\}< R< \max\{ 4|\D|^{-1}, 4(1 + 2|\D|)\} < R' /2$.

\begin{proof}
We take formula (\ref{Fprob1}) with radius $R$ as given in the conditions in the Lemma and radius $r>0$ so that $\max\{ 4|\D|^{-1}, 4(1 + 2|\D|)\}<r^{-1} < R' /2$. Note that the conditions on the contours $\cC_R$ and $\cC_r$ given in \textsc{Theorem 2} (i.e.~all the poles lie inside/outside of the contours) are satisfied for our choice of radii. Then, we deform the contour in (\ref{Fprob1}) for the $\zeta$-variables to a large radius $R'>0$, with $R'$ satisfying the conditions given in the Lemma. We begin by deforming the contour for $\zeta_N$, then the contour for $\zeta_{N-1}$, and continue successively until we deform the contour for $\zeta_1$. When we deform the contour for the $\zeta_n$ variable, we encounter three types of poles
\begin{equation}
    (a)\, 1 - \xi_i \zeta_n= 0; \quad (b)\, 1 + \zeta_i \zeta_n -2 \Delta \zeta_i = 0, \, i < n;\quad (c)\, 1 +\zeta_n \zeta_j - 2 \Delta \zeta_n, \, n < j
\end{equation}
for any $i, j =1, \dots, N$. The contribution for a type $(a)$ pole is given by the contour integral with respect to the variable $\z_n$ with contour $- \cC^{(i)}$, i.e.~a negatively oriented circle centered at $\x_i^{-1}$ with radius $r'>0$ as given in the conditions of the Lemma. Note that the only pole, with respect to the variable $\z_n$, inside the contour $-\cC^{(i)}$ is given by $\z_n = \x_i^{-1}$ because $r'$ is chosen to be small enough. The result then follows by showing that the type $(b)$ and $(c)$ poles contribute no residue.

Assume we have already deformed the $\zeta_j$ variables for $j >n$ so that $\zeta_i \in \mathcal{C}_r$ for $i <n$ and $\zeta_j \in \Omega(\xi) $ for $j >n$. We then deform the contour for the $\zeta_n$ variable. Below, we consider the residue contribution from the type $(b)$ and $(c)$ poles.

\noindent \textbf{Case (b).} We compute the residue at
\begin{equation}
    \zeta_n = (2 \Delta \zeta_{\ell} - 1)/\zeta_{\ell}
\end{equation}
for $\ell <n$. The result is a $(2N-1)$-fold contour integral with the same integrand, say $I_N(\xi, \zeta; t)$, except that the term $1 + \zeta_{\ell} \zeta_n - 2 \Delta \zeta_{\ell}$ is replaced by $\zeta_{\ell}$ and the variable $\zeta_n$ is evaluated at $(2 \Delta \zeta_{\ell} - 1)/\zeta_{\ell}$ for the rest of the terms. 

We then compute the integral with respect to the $\zeta_{\ell}$ variable for the resulting residue term. The integral is computed by analyzing the poles and residues inside the contour $\cC_{r}$ for $\zeta_{\ell}$. The possible poles are given by
\begin{equation}
    \begin{split}
    &1 - \xi_k \zeta_n = 0, \quad \hspace{11mm} k=1, \dots, N\\
    &1 + \zeta_{\ell} \zeta_j - 2 \Delta \zeta_{\ell} = 0, \quad \ell < j,\, \zeta_j \in \Omega(\xi)\\
    &1 + \zeta_i \zeta_{\ell} - 2 \Delta \z_{i}=0, \quad i < \ell, \, \z_i \in \cC_r\\
    &1 + \zeta_i \zeta_n - 2 \Delta \zeta_i =0, \quad i < n, \, \zeta_i \in \cC_r\\
    &1 + \zeta_n \zeta_j + 2\Delta \zeta_n = 0, \quad n < j, \, \zeta_j \in \Omega(\x)\\
    &\zeta_{\ell}^{x-y_j-1} \zeta_n^{x- y_n-1} =0, \quad j \neq n .
    \end{split}
\end{equation}
In particular, the location of the possible poles is given by the following
\begin{equation}
    \begin{split}
    \zeta_{\ell} = \frac{\x_k}{2 \Delta \x_k -1}  \quad &\Rightarrow \quad    \left| \frac{\x_k}{2 \Delta \x_k -1} \right| > r \\
    \zeta_{\ell} = \fr{1}{2 \Delta - \zeta_j} \quad &\Rightarrow \quad \zeta_n = 2 \Delta - \z_{\ell}^{-1} = \zeta_j\\
    \zeta_{\ell} = 2 \Delta - \z_i^{-1} \quad &\Rightarrow\quad \left| 2\Delta - \z_i^{-1}\right| >r
    \\
    \zeta_n = 2\Delta - \zeta_i^{-1} \quad &\Rightarrow \quad \zeta_{\ell} = \zeta_i\\
    \zeta_{\ell} = \fr{2\D - \z_j}{4\D^4- 2 \D \z_j -1} \quad &\Rightarrow \quad \left| \fr{2\D - \z_j}{4\D^2- 2 \D \z_j -1}\right|>r\\
    \zeta_{\ell}^{y_n - y_{\ell}} \quad &\Rightarrow \quad y_n- y_{\ell}\geq 1.
    \end{split}
\end{equation}
We use the assumptions on the radii $R, R', r'>0$ given in the statement of the Lemma and the condition on the radius $r>0$ fixed at the beginning of the proof to establish the inequalities above. For the first two inequalities, it suffices to have $R , r^{-1} > 1 + 2 |\D|$. For the third inequality, we have to consider two cases $\z_j \in \cC^{(0)}$ or $\z_j \in \cC^{(k)}$ with $k\neq 0$. In the first case when $\z_j \in \cC^{(0)}$, we have that $|\z_j| = R'$ and we use the bounds $R'>16 |\D|$ and $R' > 8 |\D|^{-1}$ that follow from the condition on the statement of the Lemma. In the second case when $\z_j\in\cC^{(k)}$ with $k \neq 0$, we have that $|\z_j| \leq (3/2)R^{-1}$ and we use the bound $(3/2)R^{-1} < |\D|$ that follows from the statement of the Lemma. Then, in all the cases above except for the second and fourth case, the poles lie outside the contour $\mathcal{C}_r$, meaning that there is no residue contribution. In the second and fourth cases, the determinant term $D_N(\xi, \zeta)$ vanishes because two columns in the matrix of the determinant are equal to each other since two $\zeta$ variables are equal to each other. In the last case, there is no pole since the exponent is positive. Then, the pole from the denominator and the zero from the determinant cancel out, meaning that these cases don't produce a residue. 

Therefore, by computing the integral with respect to the $\zeta_{\ell}$ variable, we have that the residues from the type $(b)$ poles vanish.

\noindent\textbf{Case (c).} We compute the residue at 
\begin{equation}
    \z_n = \frac{1}{2 \D - \z_{\ell}}
\end{equation}
with $n <\ell$. The result is a $(2N-1)$-fold contour integral with the same integrand, say $I_N(\xi, \zeta; t)$, except that the term $1 + \zeta_{n} \zeta_{\ell} - 2 \Delta \zeta_{n}$ is replaced by $2\D-\z_{\ell}$ and the variable $\zeta_n$ is evaluated at $1/(2\D - \z_{\ell})$ for the rest of the terms. 

In this case, we have have $\zeta_{\ell} \in \Omega(\xi)$ since $\ell >n$. Thus, we have two possibilities: (i) $\zeta_{\ell} \in \cC^{(0)}= \cC_{R'}$ , or (ii) $\zeta_{\ell} \in - \cC^{(k)}$ for some $k =1, \dots, N$ (i.e.~a negatively oriented small circle of radius $r'$ centered at $\xi_k^{-1}$). In the first case, we will not cross a pole in the contour deformation and there will be no residue to consider. In the second case, the pole will cancel out with a zero from the numerator and, again, there will be no residue to consider. We give more details below. 

In the first case, when $\zeta_{\ell} \in \mathcal{C}_{R'}$, we have $\zeta_n = 1/(2 \Delta - \zeta_{\ell})$. This pole lies inside the contour $\mathcal{C}_r$ since $R'\,r > 2 $ and $r < (1 + 2 |\D|)^{-1} $. Thus, we don't cross this pole when we deform the $\mathcal{C}_r$ contour to $\cC^{(0)} = \cC_{R'}$.

In the second case, when $\zeta_\ell  \in - \cC^{(k)}$, we first compute the residue at $\zeta_{\ell} = \xi_k^{-1}$. We obtain an $(2N-1)$-fold contour integral with the same integrand, say $I_N(\xi, \zeta; t)$, except that the determinant $D_N(\xi, \zeta)$ is replaced the same determinant with the $k^{th}$ row and the ${\ell}^{th}$ removed and multiplied by the factor $(1 + \xi_k^2 - 2\Delta \xi_k)^{-1}$, and the rest of the terms are the same with the variable $\zeta_{\ell}$ evaluated at $\xi_k^{-1}$.

We now deform the contour for $\zeta_{n}$ to the contour $\cC^{(0)}$. After taking the $\zeta_{\ell} = \xi_k^{-1}$ residue, it turns out that the terms giving rise to the pole $\zeta_n= 1/(2\Delta - \zeta_{\ell})$ becomes
\begin{equation}
    (1 +  \zeta_n \zeta_{\ell} - 2\Delta \zeta_n) = \xi_k^{-1}(\xi_k + \zeta_n - 2\Delta \xi_k \zeta_n).
\end{equation}
Note that this term also appears in the numerator of the integrand, meaning that this term cancels out and there is no residue in this case.

Therefore, when we deform the contour for the $\zeta_n$ to infinity, we don't cross any type $(c)$ poles. Moreover, this, along with the argument for the type $(b)$ poles, means that we only cross the poles due to the type $(a)$ poles. This establishes the result.

\end{proof}

\subsection{Series Expansion}

We write the contour formula (\ref{Fprob2}) as a summation by expanding the integrals over the contour $\Omega(\xi)$, given by (\ref{e:deformation_contour}), as a summations of $N+1$ integrals. We introduce some notation to encode the different terms in the summation.

Take the set of all maps from the index set $\{ 1, \dots, N\}$ to the set $\{0, 1, \dots, N\}$ and denote it by
\begin{equation}\label{e:maps}
    \mathcal{T} := \{\tau : \{1, \dots, N\} \rightarrow \{0, 1, \dots, N\} \} = \mathrm{Hom}\left(\{1, \dots, N\}, \{0,1, \dots, N\}\right).
\end{equation}
In the following, a map $\tau \in \mathcal{T}$ will correspond to a term with contours $\cC^{(\tau(k))}$, given by (\ref{e:deformation_contour}), for the $\zeta_k$ variable and $k=1, \dots, N$. Moreover, we will show that some contour integrals will vanish for certain $\tau \in \mathcal{T}$. We consider the set of maps that map injectively to the elements $\{1, \dots, N\}$ in the image and the cardinality of the preimage $\sigma^{-1}(0)$ is fixed; 
\begin{equation}\label{e:maps_cont}
    \mathcal{T}_n := \{\tau \in \mathcal{T} \mid |\tau^{-1}(0)| =n; |\tau^{-1}(k)| \leq 1,\, k =1, \dots, N  \}.
\end{equation}

\textsc{Lemma 7.2}
For $\D \neq 0$, $\cF_N(x,t)=\bP_Y(X_1(t)\ge x)$ equals
\be
\sum_{n=0}^N \sum_{\tau \in \mathcal{T}_n}\oint_{\cC_R}\cdots\oint_{\cC_R}  \oint_{\cC^{(\tau(1))}} \cdots \oint_{\cC^{(\tau(N))}}I_N(\xi, \zeta; x, t)\, d^N\z \,d^N\x\hspace{5ex}\label{Fprob3}\ee
where the integrand $I_N(\xi, \zeta; x, t)$ is the same integrand as in (\ref{Fprob2}), the summation is take over the set of maps $\mathcal{T}_n$ given by (\ref{e:maps_cont}), the contour $\cC_R$ is a circle centered at zero with radius $R>0$, the contours $\cC^{(\tau(k))}$ are given by (\ref{e:deformation_contour}) with radii $r' = R^{-1}/2, R' >0$ so that the radii satisfy the bounds $\max\{2 |\D|^{-1}, 2(1+2|\D|) \} < R < \max \{4 |\D|^{-1}, 4(1+2|\D|) \}< R'/2$ 

\begin{proof}
We take the contour formula (\ref{Fprob2}) from \textsc{Lemma7.1}. We then expand the integrals over the contours $\Omega(\xi)$ as a sum of $N+1$ integrals with contours given by the right side of (\ref{e:deformation_contour}). The result is a summation over the set of maps $\mathcal{T}$ given by (\ref{e:maps}),
\begin{equation}\label{e:sum1}
    \cF_N(x,t) = \sum_{\tau \in \mathcal{T}} \oint_{\cC_R}\cdots \oint_{\cC_R}  \oint_{\cC^{(\tau(1))}} \cdots \oint_{\cC^{(\tau(N))}}I_N(\xi, \zeta; x, t)\, d^N\z \,d^N\x.
\end{equation}
The result of this lemma follows by showing that some terms vanish, i.e.~if $\tau \notin \mathcal{T}_n$ the corresponding contour integral will vanish. Below, we show that a term in the summation vanishes if $\tau(j) = \tau(k) >0$ with $j \neq k$.

Take $\tau \in \mathcal{T}$ with $\tau(j') = \tau(k') =\ell >0$ with $n \neq m$ and $j', k'=1, \dots, N$. We show that the term in the summation (\ref{e:sum1}) with this $\tau \in \mathcal{T}$ vanishes by taking the integrals with respect to the variables $\zeta_{j'}$ and $\zeta_{k'}$. We take the integral with respect to the $\zeta_{j'}$ and $\zeta_{k'}$ variables by taking the residues at the poles given by $\zeta_{j'} = \xi_{\ell}^{-1}$ and $\zeta_{k'} = \xi_{\ell}^{-1}$. Note that the poles given by $\zeta_{j'} = \xi_{\ell}^{-1}$ and $\zeta_{k'} = \xi_{\ell}^{-1}$ correspond to the $(\ell, j')$-entry and the $(\ell, k')$-entry of the matrix for the $D_N(\xi, \zeta)$ determinant. First, we take the residue at $\zeta_{j'} = \xi_\ell^{-1}$, the determinant transforms as follows
\begin{equation}\label{e:transform1}
    D_N(\xi, \zeta) = \det \left( \frac{1}{(1 - \xi_j \zeta_k)(\xi_j+ \zeta_k - 2 \Delta \xi_j \zeta_k)} \right)_{j, k=1}^N \longrightarrow \frac{(-1)^{\tau(j')-j'-1}}{1+\xi_{\ell}^2 - 2\Delta \xi_{\ell}} \det \left( \frac{1}{(1 - \xi_j \zeta_k)(\xi_j+ \zeta_k - 2 \Delta \xi_j \zeta_k)} \right)_{j \neq \ell, k \neq j'}.
\end{equation}
For the rest of the factors in the integrand, one evaluates $\zeta_{j'} = \xi_\ell^{-1}$ when we take the residue at $\zeta_{j'} = \xi_{\ell}^{-1}$. One may check that this doesn't introduce any poles with respect to the $\zeta_{j'}$ variable inside the $\cC^{(\ell)}$ contour. Then, the residue at $\zeta_{j'} = \xi_{\ell}^{-1}$ doesn't have a pole at $\zeta_{k'} = \xi_{\ell}^{-1}$ since the pole at $\zeta_{k'} = \xi_{\ell}^{-1}$ is removed when we take the residue and no other pole is introduced. Thus, by taking the residue at $\zeta_{k'} = \xi_{\ell}^{-1}$ after taking the residue at $\zeta_{j'} = \xi_{\ell}^{-1}$, we have that the term vanishes. That is,
\begin{equation}\label{e:vanish_term1}
    \oint_{\cC_R}\cdots \oint_{\cC_R}  \oint_{\cC^{(\tau(1))}} \cdots \oint_{\cC^{(\tau(N))}}I_N(\xi, \zeta; x, t)\, d^N\z \,d^N\x = 0
\end{equation}
if $\tau(n) = \tau(m) =\ell >0$ with $n \neq m$ and $n, m=1, \dots, N$.

The result of the lemma then follows by taking the summation representation given by (\ref{e:sum1}) and noting that the terms with $\tau \notin \mathcal{T}_n$ for some $n =0, 1, \dots, N$ vanish due to the identities (\ref{e:vanish_term1}).
\end{proof}

\subsection{Residue Computations}

We compute the contour integrals with respect to the $\z_k$ variables with $\tau(k)\neq 0$ for each of the terms in the series expansion given by (\ref{Fprob3}). First, we introduce some notation to represent the resulting integrand after the residue computations.

Fix $\tau \in \mathcal{T}_{N-M}$ with $0 \leq M \leq N$ and $\mathcal{T}_{N-M}$ given by (\ref{e:maps_cont}). Then, define the following sets
\begin{equation}\label{e:index_sets}
    K_1 := \tau^{-1}(0), \quad K_2 := K_1^{c} = \{k_1 < \cdots < k_M \}, \quad J_2=\tau \left( K_2 \right) = \{\tau_1 = \tau(k_1), \dots, \tau_{M}=\tau(k_M) \}, \quad J_1 = J_2^c.
\end{equation}
We also introduce the following functions
\begin{equation}\label{e:integrands}
    \begin{split}
    I_N(\xi, \zeta; \tau)&=\frac{\prod_{j\in J_1, k \in K_1}(\x_j + \z_k - 2\D \x_j \z_k)D_N(\x, \z; \tau)}{\prod_{\substack{j<k\\ j, k \in J_1}}(1 + \x_j \x_k - 2\D \x_j) \prod_{\substack{j<k\\j,k \in K_1}}(1 +\z_j \z_k - 2\D \z_j)} \prod_{j \in J_1}\x_j^{x-y_j-1} e^{-\mi t \epsilon(\x_j)} \prod_{k \in K_1}\z_k^{x-y_k-1} e^{\mi t \epsilon(\z_k)}\\
    f(\x, \z; \tau)&= \prod_{\ell=1}^M\left(\prod_{\substack{\tau_{\ell} < k \\ k \neq \tau_{\ell+1}, \dots, \tau_M}}\left(\fr{1 + \x_{\tau_{\ell}} \x_k - 2\D \x_k}{1 + \x_{\tau_{\ell}} \x_k - 2\D \x_{\tau_{\ell}}}\right)\prod_{\substack{k_{\ell} < k \\ k \neq k_{\ell+1}, \dots, k_M}}\left(\fr{\x_{\tau_{\ell}} +  \z_k - 2\D \x_{\tau_{\ell}} \z_k}{\x_{\tau_{\ell}} +  \z_k - 2\D }\right) \right) \prod_{\ell=1}^{M}\x_{\tau_{\ell}}^{y_{k_{\ell}}- y_{\tau_{\ell}}-1}\\
    D_N(\x, \z; \tau) &= (-1)^{\sum_{\ell=1}^M \tau_{\ell} - k_{\ell}} \det\left(d(\x_j, \z_k)\right)_{j \in J_1, k \in K_1} = (-1)^{\sum_{\ell=1}^M \tau_{\ell} - k_{\ell}} \sum_{\gamma: K_1 \rightarrow J_1} \prod_{k \in K_1}(-1)^{\gamma(k)- k} d(\x_{\gamma(k)}, \z_{k})
    \end{split},
\end{equation}
where the function $d(\x, \z)$ is given by (\ref{IK}), the sum on the last line is taken over all bijections $\gamma: K_1 \rightarrow J_1$, and the sets $K_1$, $J_1$ are given by (\ref{e:index_sets}) and $M = N - |\tau^{-1}(0)|$. Note that $I_N(\x, \z; \tau)$ is equal to the integrand of contour integrals (\ref{Fprob1}), (\ref{Fprob2}) and (\ref{Fprob3}) if $|\tau^{-1}(0)| = N$.

\textsc{Lemma 7.3} Fix $\tau \in \mathcal{T}_{N-M}$, with $0 \leq M \leq N$, and take the notation from (\ref{e:index_sets}). Then, for $\D \neq 0$, we have
\begin{equation}
    \oint_{\cC_R}\cdots \oint_{\cC_R} \oint_{\cC^{\tau(1)}} \cdots \oint_{\cC^{\tau(N)}} I_N(\x, \z) d^N \z d^N \x = \oint_{\cC_R}\cdots \oint_{\cC_R} \oint_{\cC_{R'}} \cdots \oint_{\cC_{R'}} I_N(\x, \z; \tau) f(\x, \z; \tau) \left(\prod_{k \in K_1} d\z_j \right) d^N \x
\end{equation}
where the integral on the left side is a $2N$-fold contour intergal and the integral on the right side is a $(N+|K_1|)$-fold contour integral, the integrand on the left side equal to the integrand in (\ref{Fprob2}) and the integrand on the right side is given by (\ref{e:integrands}), and the contours are the same as in the statement of \textsc{Lemma 7.2} so that the $\z$ variables are integrated with respect to $\cC_{R'}$ contours.

\begin{proof}
We obtain the identity in this lemma by computing the integrals with respect to the $\z_{k_{\ell}}$ variables with $k_{\ell} \in K_2$. In particular, the contours are given by $-\cC^{(\tau_{\ell})}$, which are negatively oriented circles of radius $r'= 1/(2 R)$ and centered at $\x_{\tau_{\ell}}^{-1}$, for the integrals with respect to $\z_{k_{\ell}}$ and $k_{\ell} \in K_2$. Then, we compute the integrals by taking the residues at $\z_{k_{\ell}} = \x_{\tau_{\ell}}^{-1}$. We start by taking the residue at $\z_{k_{M}} = \x_{\tau_{M}}^{-1}$ and continue successively until we take the residue at $\z_{k_1} = \x_{\tau_1}^{-1}$.

Let's take the residue with respect to $\z_{k_M} = \x_{\tau_M}^{-1}$. Note that the pole corresponding to this residue comes from the $(\tau_M, k_M)$-entry of the matrix of the $D_{N}(\x, \z)$ determinant. Then, when we take the residue, the determinant is replaced by a determinant of the same matrix with the $\tau_M$-row and $k_M$-column removed and a prefactor $(-1)^{\tau_M-k_M}(1 + \x_{\tau_M}^2 - 2\D\, \x_{\tau_M})^{-1}$. That is,
\begin{equation}
    \det \left( \frac{1}{(1 - \xi_i \zeta_j)(\xi_i+ \zeta_j - 2 \Delta \xi_i \zeta_j)} \right)_{i, j=1}^N \longrightarrow \frac{(-1)^{\tau_M-k_M-1}}{1+\xi_{\tau_M}^2 - 2\Delta \xi_{\tau_M}} \det \left( \frac{1}{(1 - \xi_i \zeta_j)(\xi_i+ \zeta_j - 2 \Delta \xi_i \zeta_j)} \right)_{i \neq \tau_M, j \neq k_M}.
\end{equation}
The other terms of the integrand, when we compute the residue, transform by evaluating $\z_{k_M} = \x_{\tau_M}^{-1}$. Then, the result after taking the residue is
\begin{equation}
    \begin{split}
    &\frac{\prod_{j \neq \tau_M, k \neq k_M}(\x_j + \z_k - 2\D \x_j \z_k)}{\prod_{\substack{j<k\\ j, k \neq \tau_M}}(1 + \x_j \x_k - 2\D \x_j) \prod_{\substack{j<k\\j,k \neq k_M}}(1 +\z_j \z_k - 2\D \z_j)} \prod_{j \neq \tau_M}\x_j^{x-y_j-1} e^{-\mi t \epsilon(\x_j)} \prod_{k \neq k_M}\z_k^{x-y_k-1} e^{\mi t \epsilon(\z_k)}\\
    \times & (-1)^{\tau_M- k_M}\det\left(\fr{1}{(1- \x_j \z_k)(\x_j +\z_k -2\D \x_j \z_k)}\right)_{j \neq \tau_M, k \neq k_M}\\
    \times&\x_{\tau_{M}}^{y_{k_{M}}- y_{\tau_{M}}-1} \prod_{\tau_{M} < k }\left(\fr{1 + \x_{\tau_{\ell}} \x_k - 2\D \x_k}{1 + \x_{\tau_{\ell}} \x_k - 2\D \x_{\ell}}\right)\prod_{k_{M} < k }\left(\fr{\x +  \z_k - 2\D \x_{\ell} \z_k}{\x_{\ell} +  \z_k - 2\D }\right)
    \end{split}.
\end{equation}
The sign infront of the determinant changed by negative one since we are taking the integral over a negatively oriented circle.

We continue taking the integrals with respect to the variables $\z_{k_{\ell}}$, successively with $\ell$ decreasing, and evaluating the residues at $\z_{k_{\ell}} = \x_{\tau_{\ell}}^{-1}$. The computations are similar to the base case $\z_{k_{M}} = \x_{\tau_{M}}^{-1}$. In particular, the pole giving rise to residue comes from the $(\tau_{\ell}, k_{\ell})$-entry of the determinant. Then, when we take the residue, the determinant transforms by removing the $\tau_{\ell}$-row and the $k_{\ell}$-column and adding a prefactor. The other terms in the integrand transform by evaluating $\z_{k_{\ell}} = \x_{\tau_{\ell}}^{-1}$. We skip the details here since the computations are very similar to the base case. The result follows by computing all the integrals with respect to the $\z_{k_{\ell}}$ variable with $k_{\ell} \in K$.
\end{proof}

\textsc{Theorem 4.} For $\D \neq 0$, $\cF_N(x,t)=\bP_Y(X_1(t)\ge x)$ equals
\be
\sum_{n=0}^N \sum_{\tau \in \mathcal{T}_n}\oint_{\cC_R}\cdots \oint_{\cC_R} \oint_{\cC_{R'}} \cdots \oint_{\cC_{R'}} I_N(\x, \z; \tau) f(\x, \z; \tau)\, \left(\prod_{k \in K_1} d\z_k \right)\, d^N \x\hspace{5ex}\label{Fprob4}
\ee
where the integrand is given by (\ref{e:integrands}), the summation is take over the set of maps $\mathcal{T}_n$ given by (\ref{e:maps_cont}), and the contours $\cC_R$ and $\cC_{R'}$ are circles centered at zero with radii $R, R'>0$ so that $\max\{2 |\D|^{-1}, 2(1 + 2|\D|)\} < R < \max\{4 |\D|^{-1}, 4(1 + 2|\D|)\} < R'/2$.

\begin{proof}
The result is a direct consequence of \textsc{Lemma 7.2} and \textsc{Lemma 7.3}.
\end{proof}

\subsection{Deformation to Steep Descent Contours}

We take the series expansion formula (\ref{Fprob4}) and deform the contours to the steep descent contours given by (\ref{e:steep_contours}).

\begin{figure}[h!]
    \centering
    \includegraphics[scale=0.35]{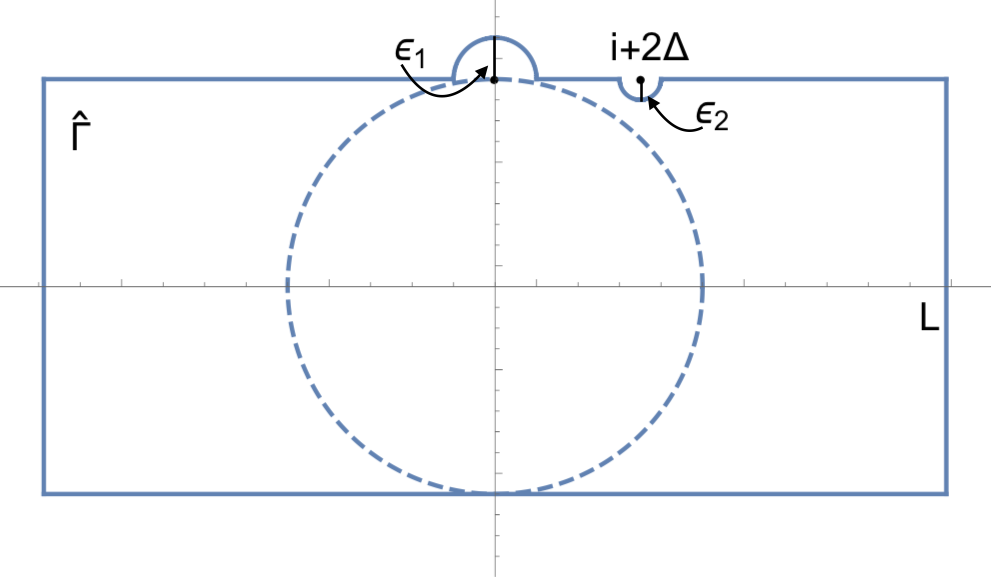}
    \caption{The contour $\widehat{\Gamma}$.}
    \label{f:rec_contour}
\end{figure}

Let $\widehat{\Gamma}$ be a positive oriented rectangle centered at zero, with length equal to $2 L = 2 \sqrt{R^2-1}$ and height equal to $2$, and two half-circle bumps as indicated on Figure \ref{f:rec_contour}. The bump centered at $\mi$ has radius $\epsilon_1$ and the bump centered at $\mi + 2\D$ has radius $\epsilon_2$ so that $0 < \epsilon_2 \ll \epsilon_1 \ll 1$. 

\textsc{Lemma 7.5} Fix $\tau \in \mathcal{T}_{N-M}$, with $0 \leq M \leq N$, take the notation from (\ref{e:index_sets}). Then, for $\D \neq 0$, we have
\begin{equation}
    \begin{split}
    &\oint_{\cC_R}\cdots \oint_{\cC_R} \oint_{\cC_{R'}} \cdots \oint_{\cC_{R'}} I_N(\x, \z; \tau) f(\x, \z; \tau)\, d^J \z\, d^N \x= \oint_{\Gamma_{+}} \cdots \oint_{\Gamma_{-}} I_{N}(\xi, \zeta; \tau)\, \left(  \oint_{\widehat{\Gamma}} \cdots \oint_{\widehat{\Gamma}}f(\x, \z; \tau) d^{J_2} \x\right) d^{K_1}\z\, d^{J_1}\x
    \end{split}
\end{equation}
where the integrand is given by (\ref{e:integrands}), the differentials $d^S \x$ or $d^S \z$ are $|S|$-fold differential over the variables $\x_s$ or $\z_s$ with $s\in S$, the contours $\Gamma_{\pm}$ are given by (\ref{e:steep_contours}) with $R_{+} = R$ and $R_{-} = R'$ so that $\x_{j} \in \Gamma_{+}$ and $\z_{k} \in \Gamma_{-}$ for $j\in J_1$ and $k \in K_1$, the contour $\widehat{\Gamma}$ is given by Figure \ref{f:rec_contour}, the contours $\cC_R$ and $\cC_{R'}$ are circles centered at zero with radii $R,R' >0$ so that $\max\{2|\D|^{-1}, 2(1 + 2|\D|)\} < R < \max\{4|\D|^{-1}, 4(1+2 |\D|)\}< R'/2$.

\begin{proof}
We obtain the result by deforming the contours and showing that we don't cross any poles. We begin by deforming the contour, for the $\x_j$ variables with $j \in J_2$, from $\cC_R$ to $\widehat{\Gamma}$. Then, for the $\x_j$ variables with $j \in J_1$, we deform the contours $\cC_R$ to the contours $\Gamma_{+}$. Finally, for the $\z_k$ variables, we deform the contours $\cC_{R'}$ to the contours $\Gamma_{-}$.

Consider the integral with respect to $\xi_{\ell} \in \cC_R$ and $\ell \in J_2$. We deform the contour $\cC_R$ to the contour $\widehat{\Gamma}$. Note that the factor $I_N(\x,\z; \tau)$ is independent of the $\x_{\ell}$ variable. Then, the only possible poles are given by
\begin{equation}\label{e:poles1}
    1 + \x_{\ell} \x_k -2\D \x_{\ell}= 0 , \quad \x_{\ell} + \z_k - 2\D =0.
\end{equation}

In the first case of (\ref{e:poles1}), the location of the pole is given by $(2\D - \x_k)^{-1}$ with $\x_k \in \cC_R$ or $\x_k \in \widehat{\Gamma}$, depending on the index and if the contour for the variable has been deformed. If $\x_k \in \cC_R$, the location of the pole $(2\D - \x_k)^{-1}$ clearly lies inside the unit circle since $R > 1+ 2|\D|$. In particular, we don't cross this pole when we deform from the contour $\cC_R$ to the contour $\widehat{\Gamma}$, since the contour $\widehat{\Gamma}$ lies outside the unit circle. If $\x_k \in \widehat{\Gamma}$, we note that the location of the pole $(2 \D - \x_k)^{-1}$ also lies inside the unit circle except for the region with the small half-circle bump of radius $\epsilon_2$. We then consider $\x_k$ lying on the small half-circle bump of $\widehat{\Gamma}$ and we write $\x_k = \mi + 2 \D + \epsilon_2\, e^{\mi \phi}$. Then, the location of the pole is given by
\begin{equation}
    (2 \D - \x_k)^{-1} = (- \mi - \epsilon_2\, e^{\mi \phi})^{-1} = \mi - \epsilon_2\, e^{\mi \phi} + \mathcal{O}(\epsilon_2^2),
\end{equation}
where the last equality follows from $0 < \epsilon_2 \ll 1$. Moreover, since $\epsilon_2 \ll \epsilon_1$, we have that the location of the pole $(2 \D - \x_k)^{-1}$ lies inside the large bump of the contour $\widehat{\Gamma}$, when $\x_k$ lies on the small bump. Then, we have that the pole $(2 \D - \x_k)^{-1}$ lies inside the unit circle if $\x_k$ doesn't lie on the small bump, and the pole lies inside the large bump if $\x_k$ lies on the small bump. In particular, if $\xi_k \in \cC_R \cup \widehat{\Gamma}$, the location of the pole lies inside the contour $\widehat{\Gamma}$ and we don't cross any poles, given by the first case of (\ref{e:poles1}), when we deform form the contour $\cC_R$ to the contour $\widehat{\Gamma}$.

In the second case of (\ref{e:poles1}), the location of the pole is given by $2 \D - \z_k$. Additionally, we have that $\z_k \in \cC_{R'}$. Given the conditions on the radii $R, R'>0$, it follows that $R < R' -2 |\D|$. Then, the pole given by $2 \D - \z_k$ lies outside the contour $\cC_R$. In particular, we don't cross the pole when we deform the contour form $\cC_R$ to $\widehat{\Gamma}$. Thus, we don't cross any poles, given by the second case of (\ref{e:poles1}), when we deform the contours from $\cC_R$ to $\widehat{\Gamma}$.

Consider now the integral with respect to $\x_{\ell} \in \cC_{R}$ with $\ell \in J_1$. We deform the contour $\cC_{R}$ to the contour $\Gamma_{+}$ with $\x_j \in \widehat{\Gamma}$ for $j \in J_2$. The location of the possible poles are given by
\begin{equation}\label{e:poles3}
    (2\D - \x_j)^{-1}, \quad 2\D - \x_j^{-1}, \quad \z_k^{-1}.
\end{equation}

In the first case of (\ref{e:poles3}), the variable $\x_j$ may lie on the the contours $\Gamma_{+}$ or $\cC_R$, depending on the index. In particular, if $\x_j \in \widehat{\Gamma}$, then $j = \tau_{k}$ for some $k$, see (\ref{e:index_sets}). Moreover, $I_N(\x, \z; \tau)$ is independent of $\x_{j} = \x_{\tau_{k}}$ and the pole due to the $f(\x, \z; \tau)$ function is of the form $(2\D- \x_{\tau_{k}})^{-1}$; see (\ref{e:integrands}). Thus, for the first case, $\x_j$ will never lie on the contour $\widehat{\Gamma}$ and only lie on the contours $\cC_{R}$ or $\Gamma_{+}$. If $\xi_j \in \cC_R$, the location of the pole $(2\D - \x_j)^{-1}$ clearly lies inside the unit circle since $R- 2|\D| >1$. In particular, we don't cross this pole when we deform from the contour $\cC_R$ to the contour $\Gamma_{+}$, since the contour $\Gamma_{+}$ lies outside the unit circle. If $\xi_j \in \Gamma_{+}$, the location of the pole $(2\D - \x_j)^{-1}$ will also lie outside the unit circle. This due to the fact the $\D$ is a real number and $R-2|\D| >1$. In particular, if $\x_j \in \cC_R \cup \Gamma_{+} $, we don't cross a pole, given by the first case of (\ref{e:poles3}) when we deform from the contour $\cC_R$ to the contour $\Gamma_{+}$.

In the second case of (\ref{e:poles3}), the variable $\x_j$ may lie on the the contours $\Gamma_{+}$, $\cC_R$, or $\widehat{\Gamma}$, depending on the index. In all three cases, we have that the $- \x_j^{-1}$ point lies inside the unit circle since the contours lie outside the unit circle. Then, the pole $2\D - \x_j^{-1}$ will lie inside $\Gamma_{+}$ since $\D$ is a real number and $2(1 + 2|\D|) < R$. In particular, if $\xi_j \in \cC_R \cup \Gamma_{+} \cup \widehat{\Gamma}$, we don't cross a pole, given by the second case of (\ref{e:poles3}), when we deform from the contour $\cC_R$ to the contour $\Gamma_{+}$.

In the third case of (\ref{e:poles3}), we have $\z_k \in \cC_{R'}$. Then, the location of the pole $\z_k^{-1}$ lies completely inside the unit circle. Then, since $\Gamma_{+}$ lies outside the unit circle, we don't cross a pole when we deform the contour $\cC_R$ to the contour $\Gamma_{+}$.

Lastly, consider the integral with respect to $\z_{\ell} \in \cC_{R'}$ with $\ell \in K_1$. We deform the contour $\cC_{R'}$ to the contour $\Gamma_{-}$. The location of the possible poles is given by
\begin{equation}\label{e:poles2}
    (2\D - \z_j)^{-1}, \quad 2\D - \z_j^{-1}, \quad 2\D- \x_k, \quad \x_k^{-1}
\end{equation}
where the variables may lie on different contours depending on the indexes.

In the first case of (\ref{e:poles2}), the variable $\z_j$ may lie on the contour $\cC_{R'}$ or on the contour $\Gamma_{-}$. In either case, the location of the pole lies completely inside the unit circle. When $\z_j \in \cC_{R'}$, this follows from the bound $R>2(1 +2 |\D|)$. When $\z_j \in \Gamma_{-}$, in addition the bound $R>2(1 +2 |\D|)$, we also need the fact that $\D$ is a real number, which means that $(2\D - \z_j)$ lies outside the unit circle for $\z_j \in \Gamma_{-}$. Then, we have that the location of the pole $(2 \D - \z_j)^{-1}$ lies completely inside the unit circle and we don't cross any poles when we deform the contour $\cC_{R'}$ to the contour $\Gamma_{-}$.

In the second case of (\ref{e:poles2}), the variable $\z_j$ may lie on the contour $\cC_{R'}$ or on the contour $\Gamma_{-}$. In either case, we know that $\z_j^{-1}$ lies inside the unit circle since $\cC_{R'}$ and $\Gamma_{-}$ lie outside the unit circle. Then, since $\D$ is a real number and $4(1+ 2|\D|) < R'$ , we have that the location of the pole $2\D -\z_j^{-1}$ lies completely inside the contour $\Gamma_{-}$. Thus, we don't cross any poles when we deform the contour $\cC_{R'}$ to the contour $\Gamma_{-}$.

In the third case of (\ref{e:poles2}), the variable $\x_k$ may lie on $\widehat{\Gamma}$ since this pole is due to the $f(\x, \z; \tau)$ factor in the integrand; see (\ref{e:integrands}). In this case, the location of the pole $2\D - \x_k$ lies completely inside the contour $\Gamma_{-}$ due to the bumps of the contour $\widehat{\Gamma}$. Since $0 < \epsilon \ll 1$, the large bump of the contour $\widehat{\Gamma}$ lies completely above the horizontal section of the contour $\Gamma_{-}$. Since the small bump in the contour $\widehat{\Gamma}$ lies inside the rectangle, the small bump will also lie completely above the V-section of the $\Gamma_{-}$ contour. Additionally, since $R + 2|\D|< R'$, the rest of the contour $\widehat{\Gamma}$ will lie completely inside the contour $\Gamma_{-}$. Then, we don't cross any poles when we deform the contour $\cC_{R'}$ to the contour $\Gamma_{-}$.

In the fourth case of (\ref{e:poles2}), we have $\x_k \in \Gamma_{+}$. Then, the location of the pole $\x_k^{-1}$ lies completely inside the unit circle, since the contour $\Gamma_{+}$ lies outside the unit circle. Then, since $\Gamma_{-}$ lies outside the unit circle, we don't cross a pole when we deform the contour $\cC_{R'}$ to the contour $\Gamma_{-}$.

We have now shown that we don't cross any poles in any case when we deform the contours. Thus, the result follows.

\end{proof}

\textsc{Proposition 7.6} For $\D \neq 0$, $\cF_N(x,t)=\bP_Y(X_1(t)\ge x)$ equals
\be
\sum_{n=0}^N \sum_{\tau \in \mathcal{T}_n}\oint_{\Gamma_{+}} \cdots \oint_{\Gamma_{-}} I_{N}(\xi, \zeta; \tau)\,\left(  \oint_{\widehat{\Gamma}} \cdots \oint_{\widehat{\Gamma}}f(\x, \z; \tau) d^{J_2} \x \right) d^{K_1}\z\, d^{J_1}\x\,\hspace{5ex}\label{Fprob5}
\ee
where the integrand is given by (\ref{e:integrands}), the sets $J_1, J_2, K_1, K_2$ are given by (\ref{e:index_sets}), the summation is take over the set of maps $\mathcal{T}_n$ given by (\ref{e:maps_cont}), and the contours $\Gamma_{\pm}$ and $\widehat{\Gamma}$ are given by (\ref{e:deformation_contour}) and Figure \ref{f:rec_contour} with $R_{+}=R, R_{-}= R'$ so that $\max\{2 |\D|^{-1}, 2(1 + 2|\D|) \} < R < \max\{4 |\D|^{-1}, 4(1 + 2|\D|) \} < R'/2$.

\begin{proof}
The result is a direct consequence of \textsc{Proposition 7.4} and \textsc{Lemma 7.5}.
\end{proof}

\section{Asymptotic Analysis, a Conjecture}\label{s:conjecture}

We believe that the formula for the probability of the left-most particle given by (\ref{Fprob4}) in \textsc{Theorem 7.6} may be suitable for asymptotic analysis when $t \ll N \rightarrow \infty$. Note that we have decomposed the integrand into two factors, $I_N(\x, \z; \tau)$ and $f(\x, \z; \tau)$. In particular, note that that the factor $f(\x, \z; \tau)$ is independent of time $t$. Additionally, for the variables of the term $I_N(\x,\z; \tau)$, we have deformed the contours to steepest descent paths. Thus, in the asymptotic limit, we expect the main contribution for the $I_N(\x,\z;\tau)$ term to come from the saddle point $(\x_0,\z_0) =(i,-i)$. Moreover, we expect the asymptotic limit of $I_N(\x, \z ; \tau)$ to be given by the Airy kernel.  We give some details of the computation below but, unfortunately, we don't give all the technical details here. The arguments below need more careful consideration.

Fix $\tau \in \mathcal{T}_n$ and let's consider the contribution of the contour integrals near the saddle point. We use the following notation for the index sets:
\begin{equation}
\begin{split}
    K_1 := \tau^{-1}(0) , \quad K_2 := (K_1)^c, \quad 
    J_1 := \tau(K_2)^c , \quad J_2 := \tau(K_2)
\end{split}
\end{equation}
The sets $K_1$ and $K_2$ will be used to index the $\z$-variables and the sets $J_1$ and $J_2$ will be used to index the $\x$-variables. In particular, variables with index from the sets $K_1$ and $J_1$ will lie on the contours $\Gamma_{\pm}$, respectively, and the variables with index from the set $J_2$ will lie on the contour $\widehat{\Gamma}$. There are no variables with index from the set $K_2$ because these variable have been integrated out, but nonetheless, this index set will appear in our formulas. Note $K_1 \cup K_2 = J_1 \cup J_2 = \{1, \dots, N \}$.

Recall that the spectral function $G$ and $H$, given in (\ref{e:spectral_function}), have a double critical point at $\x = i$ and $\z=-i$, respectively, when $x =-2t$. Let $\cB(z,r)$ be an open ball centered at $z \in \mathbb{C}$ of radius $r>0$ and $\cB(z,r)^c$ be its complement. Then, we take the following scaling
\begin{equation}\label{e:scaling}
    x= -2t - st^{1/3}, \quad \x = \mi + \mi\, \tilde{\xi}\, t^{-1/3}, \quad \z = -\mi + \mi\, \tilde{\z}\, t^{-1/3},\quad y_j+1 = v_j \, t^{1/3}
\end{equation}
if $\x \in \cB(i, t^{- \alpha})$ and $\z \in \cB(-i, t^{-\alpha})$ with $1/3 < \alpha < 1/4$.

We also have that the integrand $I_N(\x, \z; \tau)$ is exponentially small if $\x_j \in \cB(i, t^{-\alpha})^c$, for $j \in J_1$, or $\z_j \in \cB(-i, t^{- \alpha})^c$, for $j \in K_1$. This follows from \textsc{Lemma 6.1}. Additionally, we may uniformly bound the factor $f(\x, \z; \tau)$, independently of $t$, on all the $\x$ and $\z$ variables. Then, we may restrict the contours $\Gamma_{\pm}$ to the a neighborhood around the saddle points and only lose an exponentially small term. That is,
\begin{equation}\label{e:approx1}
    \begin{split}
    &\oint_{\Gamma_{+}} \cdots \oint_{\Gamma_{-}} I_{N}(\xi, \zeta; \tau)\, \left( \oint_{\widehat{\Gamma}} \cdots \oint_{\widehat{\Gamma}}f(\x, \z; \tau) d^{J_2} \x \right) d^{K_1}\z\, d^{J_1}\x\\
    &= \oint_{\Gamma_{+} \cap \cB(\mi, t^{-\alpha})} \cdots \oint_{\Gamma_{-} \cap \cB(-\mi , t^{-\alpha})} I_{N}(\xi, \zeta; \tau)\, \left( \oint_{\widehat{\Gamma}} \cdots \oint_{\widehat{\Gamma}}f(\x, \z; \tau) d^{J_2} \x\right) d^{K_1}\z\, d^{J_1}\x\, + \mathcal{O}(e^{- C t^{1-3\alpha}})
    \end{split}
\end{equation}
for some positive constant $C>0$, based on \textsc{Lemma 6.1}, and $1/3 < \alpha <1/4$.

Let us now approximate the integrands $I_N(\x, \z;\tau)$ and $f(\x, \z;\tau)$ when $\x_j \in \Gamma_{+} \cap \cB(i, t^{- \alpha})$, for $j \in J_1$, and $\z_k \in \Gamma_{-} \cap \cB(-i, t^{-\alpha})$, for $k \in K_1$. In particular, we take the scaling (\ref{e:scaling}) for the variables with indexes in the sets $J_1$ and $K_1$, for the $\x$-variables and $\z$-variables respectively. 

Note that $I_N(\x, \z;\tau)$ only depends on the variables with indexes from the sets $K_1$ and $J_1$. Then, we have
\begin{equation}\label{e:approx2}
    \begin{split}
    I_N(\x, \z; \tau)&= (-1)^{|\cX_1|+\sum_{j \in \cZ_2}\tau(j)-j}\sum_{\gamma : \cZ_1 \rightarrow \cX_1} \prod_{k \in \cZ_1}(-1)^{\gamma(k) - k} (\mi)^{y_{\gamma(k)}- y_k} g(\tx_{\gamma(k)}, \tz_k; v_{\gamma(k)}, v_k)\, t^{n/3} + \mathcal{O}(t^{(n-1)/3})\\
    g(\x, \z; x, z)&=  \fr{\exp\left( \fr{1}{3}\x^3 - \fr{1}{3}\z^3 - (s + x)\,\x +(s + z)\,\z\right)}{(\x - \z)},
    \end{split}
\end{equation}
where the sum is taken over all bijections $\gamma : K_1\rightarrow J_1$. This approximation is obtained by expanding the determinant in the term $I_{N}(\x, \z;\tau)$, given by (\ref{e:integrands}) and taking the scaling (\ref{e:scaling}). More details regarding this approximation are given in Appendix \ref{s:I_approx}.

Now, consider the approximation of the term $f(\x, \z;\tau)$ when $\x_j \in \Gamma_{+} \cap \cB(i, t^{- \alpha})$, for $j \in J_1$, and $\z_k \in \Gamma_{-} \cap \cB(-i, t^{-\alpha})$, for $k \in K_1$. We introduce the following function
\begin{equation}\label{e:B_fun}
    B(\x;\tau) = \prod_{\substack{j < k, j,k\in K_2 \\ \tau(k) < \tau(j)}} \left(\fr{1 + \x_{\tau(k)} \x_{\tau(j)}- 2\D \x_{\tau(j)}}{1 + \x_{\tau(k)} \x_{\tau(j)}- 2\D \x_{\tau(k)}} \right),
\end{equation}
with the indexes $j,k \in K_2$ and $\tau(k), \tau(j) \in J_2$. Also, let us denote the number of inversions of the $\tau$ map as follows,
\begin{equation}\label{e:nu1}
    \begin{split}
    \nu_1(j; \tau) &: = \# \{j' \in K_2 \mid j' < j , \quad \tau(j')>\tau(j)\ \}\\
    \nu_2(j; \tau) &: = \# \{j' \in K_2 \mid j < j' , \quad \tau(j)>\tau(j')\ \}\\
    \nu(j;\tau) &:= j-\tau(j) +\nu_2(j; \tau) - \nu_1(j; \tau)
    \end{split}.
\end{equation}
Note that, in the case $K_2 = \{1, 2, \dots, N\}$, we have $\nu(j;\tau)  = 0$ for $j=1, \dots, N$. Then, by taking the scaling (\ref{e:scaling}), we obtain
\begin{equation}\label{e:approx3}
    f(\x, \z; \tau) = B(\x;\tau) \prod_{j \in K_2}\left( \fr{\x_{\tau(j)} - (2\D +\mi)}{(2\mi \D +1) \x_{\tau(j)}- \mi}\right)^{\nu(j;\tau)} \prod_{j\in K_2} \x_{\tau(j)}^{y_{j}- y_{\tau(j)}-1} + \mathcal{O}(t^{-1/3}).
\end{equation}
This approximation is obtained by applying the scaling (\ref{e:scaling}) and taking the leading term in the $t^{-1/3}$ expansion of the $f(\x,\z; \tau)$ function. More details regarding this approximation are given in Appendix \ref{s:f_approx}. 

We now combine the approximations (\ref{e:approx1}), (\ref{e:approx2}) and (\ref{e:approx3}), given above. Note that the leading term of the approximation (\ref{e:approx3}) is independent of the $\tilde{\x}$ and $\tilde{\z}$ variables. We then introduce the term
\begin{equation}\label{e:F}
    F(\tau) = (\mi)^{|K_2|}\oint_{\widehat{\Gamma}} \cdots \oint_{\widehat{\Gamma}} B(\x;\tau) \prod_{j \in K_2}\left( \fr{\x_{\tau(j)} - (2\D +\mi)}{(2\mi \D +1) \x_{\tau(j)}- \mi}\right)^{\nu(j; \tau)} \prod_{j\in K_2} (\mi\,\x_{\tau(j)})^{y_{j}- y_{\tau(j)}-1}d^{J_2} \x,
\end{equation}
where we have taken the leading term of the $f(\x,\z;\tau)$ function and also incorporated the $(\mi)^{y_{\gamma(k)} - y_k}$ term from the approximation of $I_N(\x, \z;\tau)$ given by (\ref{e:approx2}), noting that $\sum_{k\in K_1}y_{\gamma(k)} - y_k +\sum_{k\in K_2}y_{\tau(k)} - y_k =0$. Then, for fixed $\tau \in \mathcal{T}_n$, we obtain the following approximation near the saddle point
\begin{equation}\label{e:approx4}
    \begin{split}
    &\oint_{\Gamma_{+}} \cdots \oint_{\Gamma_{-}} I_{N}(\xi, \zeta; \tau)\,\left( \oint_{\widehat{\Gamma}} \cdots \oint_{\widehat{\Gamma}}f(\x, \z; \tau) d^{J_2} \x\right) \, d^{K_1}\z\, d^{J_1}\x\\
    &=t^{-n/3}(-1)^{|J_1|+\sum_{k \in K_2}\tau(k)-k} \sum_{\gamma : K_1 \rightarrow J_1}  F(\tau) \prod_{k \in K_1}(-1)^{\gamma(k) - k} \mathbf{K}_{Ai}\left(s + v_{\gamma(k)}, s + v_k \right)\\
    &\quad \quad + \mathcal{O}(t^{(1-n)/3}) + \mathcal{O}(e^{-Ct^{1- 3 \alpha}}).
    \end{split}
\end{equation}
The $t^{-n/3}$ term and the Airy kernel $\mathbf{K}_{Ai}$ are obtained by taking the change of variables (\ref{e:scaling}) and the following expression for the Airy kernel
\begin{equation}\label{e:Airy}
    \mathbf{K}_{Ai}(x,z) = \int_{\infty\, e^{- 2\pi \mi /3}}^{\infty\, e^{2\pi \mi /3}} \int_{\infty\, e^{- \pi \mi /3}}^{\infty\,e^{ \pi \mi /3}} \fr{\exp\left(\frac{1}{3}\x^3  - \frac{1}{3}\z^3 - x\, \x + z\, \z\right)}{\x - \z}\, d\x\, d\z,
\end{equation}
where the contours for the $\x$ (resp.~$\z$) variable starts at $\infty\, e^{-\pi i /3}$ (resp.~$\infty\, e^{-2\pi i /3}$) goes through the origin and ends at $\infty\, e^{\pi i /3}$ (resp.~$\infty\, e^{2\pi i /3}$).

Let's now consider the formula (\ref{Fprob5}) and, in particular, the summation over $\mathcal{T}_n$ and $n$. We substitute the term in the summation by the right side of the approximation (\ref{e:approx4}). The result is a summation over $\mathcal{T}_n$, $n$, and injective maps $\gamma: K_1\rightarrow J_1$. More precisely, the summation is over a pair of bijective maps
\begin{equation}
    \tau: K_2 \rightarrow J_2, \quad \gamma: K_1 \rightarrow J_1,
\end{equation}
where $K_1 \cup K_2 = J_1 \cup J_2 = \{1, 2, \dots, N\}$. This means that we may write the summation, over $\mathcal{T}_n$, $n$ and the injective maps $\gamma: K_1 \rightarrow J_1$ and $\tau: K_2 \rightarrow J_2$, as the summation over permutations of the set $[N] =\{1, 2, \dots, N\}$. In particular, we may uniquely identify a pair of bijective maps $(\tau, \gamma)$ with a  permutation $\sigma \in \cS_N$ and a subset $S \subset [N]$ so that $(\tau, \gamma) = (\sigma|_{S^c}, \sigma|_{S})$, where the right side are restrictions of the permutation to the indicated sets. Then, under this identification, we rewrite some the notation introduced earlier. For $(\sigma, S)$ with $\sigma|_{S^c} = \tau$, we have
\begin{equation}
    B(\x;\sigma, S) = B(\x; \tau) = \prod_{\substack{j, k \in S^c, j<k \\ \sigma(j) > \sigma(k)}} \left(\fr{1 + \x_{\sigma(k)} \x_{\sigma(j)} - 2\D \x_{\sigma(j)}}{1 + \x_{\sigma(k)} \x_{\sigma(j)} - 2\D \x_{\sigma(k)}} \right).
\end{equation}
Additionally, for $(\sigma, S)$ with $\sigma|_{S^c} = \tau$, we write the inversion sets as follows,
\begin{equation}
    \begin{split}
    \nu_1(j; \sigma, S) &= \nu_1(j; \tau)  = \# \{j' \in S^c \mid j' < j , \quad \sigma(j')>\sigma(j)\ \}\\
    \nu_2(j; \sigma, S) &= \nu_2(j; \tau)  = \# \{j' \in S^c \mid j < j' , \quad \sigma(j)>\sigma(j')\ \}\\
    \nu(j; \sigma, S) &= \nu(j;\tau) = j-\sigma(j) +\nu_2(j; \sigma, S) - \nu_1(j; \sigma, S).
    \end{split}.
\end{equation}
Lastly, for $(\sigma, S)$ with $\sigma|_{S^c} = \tau$, we write
\begin{equation}
    F(\sigma, S) = F(\tau) = (\mi)^{|S^c|} \oint_{\widehat{\Gamma}} \cdots \oint_{\widehat{\Gamma}} B(\x;\sigma, S) \prod_{j \in S^c}\left( \fr{\x_{\sigma(j)} - (2\D +\mi)}{(2\mi \D +1) \x_{\sigma(j)}- \mi}\right)^{\nu(\sigma, S)} \prod_{j\in S^c} (\mi \, \x_{\sigma(j)})^{y_{j}- y_{\tau(j)}-1}d^{\sigma(S^c)} \x.
\end{equation}
Then, under the identification of the pair of injective maps and the permutations, we have
\begin{equation}
    \begin{split}
    &\sum_{n=0}^N \sum_{\tau \in \mathcal{T}_n}\oint_{\Gamma_{+}} \cdots \oint_{\Gamma_{-}} I_{N}(\xi, \zeta; \tau)\,\left( \oint_{\widehat{\Gamma}} \cdots \oint_{\widehat{\Gamma}}f(\x, \z; \tau) d^{J_2} \x\right)\, d^{K_1}\z\, d^{J_1}\x\\
    &=\sum_{\sigma \in \cS_N} (-1)^{\sigma} \sum_{S \subset [N]} (-1)^{|S|} t^{- |S|/3} \left(F(\sigma, S)\prod_{k \in S} \mathbf{K}_{Ai}\left( s + \fr{y_{\sigma(k)}+1}{t^{1/3}}, s + \fr{y_k+1}{t^{1/3}}\right) + \mathcal{O}(t^{-1/3}) + \mathcal{O}(e^{-C t^{1- 3\alpha}})\right)
    \end{split}.
\end{equation}

Assuming that the error terms don't contribute in the limit, we have the following conjecture.

\textsc{Conjecture 8.1} As $t\ll N \rightarrow \infty$, $\cF_N(x,t)=\bP_Y(X_1(t)\ge x)$, with $x = -2t -s \, t^{-1/3}$ and $y_j +1 = v_j\, t^{1/3}$, equals to the limit of
\begin{equation}\label{Fprob6}
    \sum_{\sigma \in \cS_N} (-1)^{\sigma} \sum_{S \subset [N]} (-1)^{|S|} t^{- |S|/3} F(\sigma, S)\prod_{k \in S} \mathbf{K}_{Ai}\left( s + v_{\sigma(k)}, s +v_k\right) 
\end{equation}
where $F$ is given by (\ref{e:F}) and the Airy kernel $\mathbf{K}_{Ai}$ is given by (\ref{e:Airy}).

At the moment, we are not able to control the limit of (\ref{Fprob6}) when $t \ll N \rightarrow \infty$. The main obstacle is the term $F(\sigma, S)$ on (\ref{Fprob6}). However, under some assumptions, we may simplify (\ref{Fprob6}) as a determinant of the difference of two kernels. For instance, assume
\begin{equation}\label{e:assump}
    F(\sigma,S) = \prod_{j \in S^c} \mathbf{Q}(\sigma(j), j)
\end{equation}
for some kernel $\mathbf{Q}$ on the set $\{1, \dots, N\}$. Then, we have
\begin{equation}
    \begin{split}
    &\sum_{\sigma \in \cS_N} (-1)^{\sigma} \sum_{S \subset [N]} (-1)^{|S|} t^{- |S|/3} F(\sigma,S)\prod_{k \in S} \mathbf{K}_{Ai}\left( s + \fr{y_{\sigma(k)}+1}{t^{1/3}}, s + \fr{y_k+1}{t^{1/3}}\right) \\
    &=\sum_{\sigma \in \cS_N} (-1)^{\sigma} \sum_{S \subset [N]} (-1)^{|S|} t^{- |S|/3} \prod_{j \in S^c} \mathbf{Q}(\sigma(j), j) \prod_{k \in S} \mathbf{K}_{Ai}\left( s + \fr{y_{\sigma(k)}+1}{t^{1/3}}, s + \fr{y_k+1}{t^{1/3}}\right)\\
    &= \sum_{\sigma \in \cS_N} (-1)^{\sigma} \prod_{k=1}^N\left(\mathbf{Q}(\sigma(k), k) - t^{-1/3}\, \mathbf{K}_{Ai}\left(s + \fr{y_{\sigma(k)}+1}{t^{1/3}}, s + \fr{y_k+1}{t^{1/3}} \right) \right)\\
    & = \det \left( \mathbf{Q}(j, k) - t^{-1/3}\mathbf{K}_{Ai}\left(s + \fr{y_{j}+1}{t^{1/3}}, s + \fr{y_k+1}{t^{1/3}} \right)\right)_{j,k=1}^N,
    \end{split}
\end{equation}
given the assumption (\ref{e:assump}). In fact, when $\Delta = 0$, one may check the assumption to be true and we have
\begin{equation}
    F(\sigma,S) = \mathds{1}\left(\sigma|_{S^c} = \mathrm{Id}_{S^c} \right) = \prod_{j \in S^c} \mathds{1}(\sigma(j) = j),
\end{equation}
where the functions with $\mathds{1}$ are indicator functions. This identity is easy to check since the first two terms in the intergand for $F(\sigma,S)$, given by (\ref{e:F}), are identically equal to one when $\Delta =0$. Then, we have
\begin{equation}
    \begin{split}
    &\det \left( \mathbf{Id}(j, k) - t^{-1/3}\mathbf{K}_{Ai}\left(s + \fr{y_{j}+1}{t^{1/3}}, s + \fr{y_k+1}{t^{1/3}} \right)\right)_{j, k=1}^N\\
    &=\sum_{\sigma \in \cS_N} (-1)^{\sigma} \sum_{S \subset [N]} (-1)^{|S|}t^{- |S|/3} F(\sigma, S)\prod_{k \in S} \mathbf{K}_{Ai}\left( s + \fr{y_{\sigma(k)}+1}{t^{1/3}}, s + \fr{y_k+1}{t^{1/3}}\right)
    \end{split}
\end{equation}
when $\Delta = 0$. This means that \textsc{Conjecture 8.1} is true when $\Delta =0$. Moreover, if $\Delta=0$ and $y_j=j$, we may take the limit $t \ll N \rightarrow \infty$. The right side becomes a sum of Riemann integrals, corresponding to the series expansion of a Fredholm determinant. Then, we have
\begin{equation}
    \begin{split}
    \lim_{N \rightarrow \infty}\bP_Y\left(\fr{X_1(t)+2t}{t^{1/3}}\ge -s\right) &= \lim_{t \ll N \rightarrow \infty}  \sum_{\sigma \in \cS_N} \sum_{S \subset [N]} t^{- |S|/3} \det \left( \mathbf{K}_{Ai}\left( s + \fr{j+1}{t^{1/3}}, s + \fr{k+1}{t^{1/3}}\right) \right)_{j, k \in S}\\
    &= \det \left(\mathbf{Id} - \mathbf{K}_{Ai} \right)_{L^2(s, \infty)}\\
    &=F_2(s).
    \end{split}
\end{equation}
This matches the earlier result (\ref{TWF2}) for $\Delta =0$.

We also may compute the terms in (\ref{Fprob6}) when $S=\emptyset$ and $S^c = \{1, \dots, N\} =[N]$. In that case, the formula for $F(\sigma, \emptyset)$ simplifies as follows
\begin{equation}
    F(\sigma, \emptyset) = \oint_{\widehat{\Gamma}} \cdots \oint_{\widehat{\Gamma}}\prod_{\substack{j < k \\ \sigma(k) < \sigma(j)}} \left(\fr{1 + \x_{\tau(k)} \x_{\tau(j)}- 2\D \x_{\tau(j)}}{1 + \x_{\tau(k)} \x_{\tau(j)}- 2\D \x_{\tau(k)}} \right)  \prod_{j=1}^N \x_{\tau(j)}^{y_{j}- y_{\tau(j)}-1}d^N \x,
\end{equation}
where $i,j= 1, 2, \dots, N$ on the first product of the integrand. Additionally, we may deform the contours $\widehat{\Gamma}$ to arbitrarily large circles centered at the origin. Note that $(-1)^{\sigma}F(\sigma,\emptyset)$ is equal to the integral inside the sum of (\ref{psi_N_large}) with $x_i =y_i$, for $i=1, \dots, N$, and $t=0$. Then, by \textsc{Theorem 1a}, we have
\begin{equation}
    \sum_{\sigma \in \cS_N} (-1)^{\sigma}F(\sigma, \emptyset) = 1
\end{equation}
for any $N>0$.

\section*{Acknowledgement}

The authors thank F.~Colomo, B.~Nachtergaele, and L.~Petrov for their helpful communications. This work was supported by the National Science Foundation under the grant DMS--1809311 (second author). This work was also supported by the National Science Foundation under Grant No. DMS--1928930 while the first author participated in the program \emph{``Universality and Integrability in Random Matrix Theory and Interacting Particle Systems"} hosted by the Mathematical Sciences Research Institute in Berkeley, California, during the Fall 2021 semester. Additionally, the first author was partially supported by the \emph{Engineering and Physical Sciences Research Council (EPSRC)} through grant EP/R024456/1.

\appendix
\section{Large Contour Formula}\label{s:appendix_a}

We give a proof for the large contour integral formula for the wave function, given in \textsc{Theorem 1a} by \eqref{psi_N_large}. The arguments in this section are almost verbatim the arguments in \cite{TW1}, Section 2. We apply the coordinate Bethe ansatz to obtain contour integral formulas for the wave function. We adapt the argument so that, instead of shrinking contours to zero, we expand contours to infinity. With this in mind, we redefine some objects like the subsets $\mathbb{S}(A)$, given by \eqref{e:s_subset}, and the change of variables given by \eqref{e:cv}. Otherwise, all else works the same.

\begin{proof}[Proof of Theorem 1a]
Denote the right side of \eqref{psi_N_large} by $u(X \rightarrow Y; t)$. We may show that the contour formula $u$ is equal to the wave function by checking the Schrodinger equation. In fact, by the Bethe ansatz, it suffices to check the delta initial conditions,
\begin{equation}
    u(X\rightarrow Y; 0) = \mathds{1}(X=Y)
\end{equation}
for the contour integral to be equal to the wave function. We check this condition in the following \textsc{Lemma A.1}.
\end{proof}

\subsection{Integral Cancellations}

The initial condition is satisfied by the summand in \eqref{psi_N_large} coming from the identity permutation $\rm  id$. So what we have to show is
\begin{equation}
    \sum_{\s\ne {\rm id}}\int_{\cC_R}\cd\int_{\cC_R} A_\s(\x) \prod_i\x_{\s(i)}^{x_i} \prod_i \left(\xi_i^{-y_{i}-1}\,
\me^{-\mi t \ve(\xi_i)}\right)
 d\x_1\cd d\x_N
\end{equation}
when $x_1 < \cdots <x_N$. We write $I(\sigma)$ for the integral corresponding to $\sigma$, so that the above becomes
\begin{equation}
\sum_{\sigma\ne{\rm id}}I(\sigma)=0.\label{Isum}
\end{equation}

For $1\le n<N$, fix $n-1$ distinct numbers $i_1,\dots,i_{n-1}\in[2,\,N]$. Define 
\[
A=\{i_1,\dots,i_{n-1}\},
\]
and then
\begin{equation}\label{e:s_subset}
\mathbb{S}_N(A)=\{\s\in\mathbb{S}_N:\sigma(N)=i_1,\dots,\sigma(N+2-n)=i_{n-1},\,\sigma(N+1-n)=1\}.
\end{equation}
When $n = 1$ this consists of all permutations with $1$ in position $N$, and when $n=N-1$ it consists of a single permutation. If $B$ is the complement of  $A\cup\{N\}$  in $[1,\,N]$,
then $\sigma \in\mathbb{S}_N(A)$ is determined by the restriction
\[
\s|_{[1,\,N-n]}:[1,\,N-n]\to B.
\] 

\textsc{Lemma A.1}. For each $A$,
\[
\sum_{\sigma \in\mathbb{S}_N(A)}I(\s)=0.
\]

\noindent{\bf Start of the Proof}. When $\sigma \in\mathbb{S}_N(A)$ the inversions involving $N$ are the $(i,1)$ with $i\in B$. Therefore the integrands in $I(\s)$ involving these $\s$ may be written

\[
\prod_{i\in B}S(\x_i,\x_1)\times \prod_{i\le N}\x_i^{x_{\si(i)}-y_i-1} \times\prod\{S(\x_k,\x_\l):k>\l>1,\ \s\inv(k)<\s\inv(\l)\}.
\]

The integrals are taken over $\cC_R$ with $R$ so large that all the denominators in the $S$-factors are nonzero on and outside the contour. In these integrals we make the substitution
\begin{equation}\label{e:cv}
\x_1\to{\fr{\e}{\prod_{1<i}\x_i}},
\end{equation}
so that $\e$ runs over a circle of radius $R^N$. The integrand becomes
\be(-1)^{N-n}\prod_{i\in B}\fr{1 + \eta^{-1} \prod_{\l \neq 1, i} \x_\l - 2 \Delta \e\inv \prod_{\l \neq 1} \x_\l}{1 + \eta^{-1} \prod_{\l \neq 1, i} \x_\l - 2 \Delta \x_i\inv}\label{first}\ee
\be\times\ \e^{x_{N+1-n}-y_1-1}\,\prod_{i<N}\x_i^{x_{\s\inv(i)}-x_{N+1-n}+y_1-y_i-1}\label{second}\ee
\be \times\ \prod\{S(\x_k,\x_\l):k>\l>1,\ \s\inv(k)<\s\inv(\l)\}.\label{third}\ee
(The reason that we still have $-1$ in the exponents in (\ref{second})
is that $d\x_1=\prod_{1<i}\x_i\inv\,d\e$.)

\noindent{\textsc Sublemma A.1}. When $n=N-1$ we have $I(\s)=0$.

\noindent{\bf Proof}. There is a single $i\in B$ and (\ref{first})  is analytic inside the $\x_i$-contour except for a simple pole at $\x_i=\infty$. The power of $\x_i$ in (\ref{second}) is 
\[\xi_i^{x_1-x_{2}+y_1-y_i-1},\]
and since $x_2>x_{1}$ and $y_i>y_{1}$, the exponent is negative. Therefore the integrand is analytic outside the $\x_i$-contour, and so the integral is zero.   

\noindent{\textsc Sublemma A.2}. When $n<N-1$ all $I(\s)$ with $\s\in\mathbb{S}_N(A)$ are sums of lower-order integrals in each of which (\ref{first}) is replaced by a factor depending on $A$. The other factors remain the same. In each integral some $\x_i$ with $i\in B$ is equal to another $\x_j$ with $j\in B$.

\noindent{\bf Proof}. We may assume that $\Delta\ne0$. This case follows by a limiting argument. We are going to expand some of the $\x_i$-contours with $i\in B$ to infinity. Due to the defining property of $R$, the only poles we pass will come from the product (\ref{first}). In fact, to avoid double poles later we take $\x_i\in \cC_{R_i}$ with the $R_i$ all slightly different.

Take $j=\min B$ and expand the $\x_j$-contour to infinity. The product (\ref{first}) has a simple pole at $\x_j= \infty$ (the $j$-factor has the pole and the $i$-factors with $i\ne j$ are analytic there) and the power of $\x_j$ in (\ref{second}) is negative as before, so the integrand is analytic at $\x_j=\infty$. For each $k\in B$ with $k\ne j$ we pass the pole at
\be\x_j=\fr{2 \Delta \x_k\inv -1}{\e\inv \prod_{\l \ne 1, j, k} \xi_\l}.\label{polej}\ee
coming from the $k$-factor in (\ref{first}). (Our assumption on the $R_i$ assures that there are no double poles.) For the residue we replace the $k$-factor by
\be -\fr{1 + \eta^{-1} \prod_{\l \neq 1, k} \x_\l - 2 \Delta \e\inv \prod_{\l \neq 1} \x_\l}{\eta^{-1} \prod_{\l \neq 1, j, k} \x_\l },\label{k}\ee
where in this and the $j$-factor we replace $\x_j$ by the right side of (\ref{polej}). When $i\ne j,k$ the $i$-factor becomes
\[\fr{1+ \e\inv\prod_{\l \neq 1, i} \x_\l - 2 \Delta \e\inv \prod_{\l \neq 1} \x_\l}{1-\x_i\inv\x_k}\]
and we replace $\x_j$ in the numerator by  the right side of (\ref{polej}).

We now expand the $\x_k$-contour. There is a pole of order 2 at 
$\x_k=\infty$ coming from (\ref{k}) and the $j$-factor in (\ref{first}). Since $1<j=\min B<k$, we have $y_1-y_k\le -2$, so the exponent of  $\x_k$ in (\ref{second}) is at most $-2$. Therefore the integrand is analytic at $\x_k=\infty$. The factor (\ref{k}) has no other poles outside $\cC_{R_k}$. An $i$-factor with $i\ne j,k$ will have a pole at $\x_k=\x_i$  if $R_i<R_k$.  There is also the pole at
\[\x_k=\fr{2 \Delta \x_i\inv -1}{\prod_{\l \ne q, k, j} \x_l}\]
coming from the $j$-factor. But this relation and (\ref{polej}) imply $\x_j=\x_k$. 

Thus when we expand the $\x_j$-contour and the $\x_k$-contours to infinity with $k\ne j$ we obtain $(N-2)$-dimensional integrals in each of which two of the $\x$-variables corresponding to indices in $B$ are equal. This proves the sublemma. 

\noindent{\textsc Sublemma A.3}. For each integral of Sublemma 2.2 there is a partition of $\mathbb{S}_N(A)$ into pairs $\s,\,\s'$ such that $I(\s)+
I(\s')=0$ for each pair.

\noindent{\bf Proof}. Consider an integral in which $\x_i=\x_j$. We pair $\s$ and $\s'$ if $\s\inv(i)={\s'}\inv(j)$ and $\s\inv(j)={\s'}\inv(i)$, and $\s\inv(k)={\s'}\inv(k)$ when $k\ne i,j$. The factor (\ref{second}) is clearly the same for both when $\x_i=\x_j$, and we shall show that the $\s$- and $\s'$-factors in (\ref{third}) are negatives of each other when $\x_i=\x_j$.

Assume for definiteness that
\be i<j\ \ {\rm and}\ \ \s\inv(i)<\s\inv(j).\label{assume}\ee
(Otherwise we reverse the roles of $\s$ and $\s'$.) Then the factor 
$S(\x_j,\x_i)$ does not appear for $\s$ in (\ref{third}) but it does appear for $\s'$.  
This factor equals $-1$ when $\x_i=\x_j$. 

To complete the proof it is enough to show that for any $k\ne i,j$ the product of $S$-factors involving $k$ and either $i$ or $j$ is the same for $\s$ and $\s'$ when $\x_i=\x_j$. There are nine cases, depending on the position of $k$ relative to $i$ and $j$ and the position of 
$\s\inv(k)$ relative to $\s\inv(i)$ and $\s\inv(j)$. 
If $k$ is outside the interval $[i,\,j]$ and $\s\inv(k)$ is outside the interval $[\s\inv(i),\,\s\inv(j)]$ then the products of $S$-factors for $\s$ and $\s'$ are clearly the same. There are five remaining cases, with the results displayed in the table below. The first column gives the position of $k$ relative to $i$ and $j$, the second column gives the position of $\s\inv(k)$ relative to 
$\s\inv(i)$ and $\s\inv(j)$, the third column gives the product of $S$-factors involving $k$ and either $i$ or $j$ for $\s$, and the fourth column gives the corresponding product for $\s'$. Keep (\ref{assume}) in mind.

\[\begin{array}{llllll}
i<k<j&\si(k)<\si(i) && S(\x_k,\x_i) && S(\x_k,\x_j)\\
i<k<j&\si(i)<\si(k)<\si(j) && 1 && S(\x_k,\x_j)\,S(\x_i,\x_k)\\
i<k<j&\si(j)<\si(k) && S(\x_j,\x_k) && S(\x_i,\x_k)\\
k<i&\si(i)<\si(k)<\si(j) && S(\x_i,\x_k) && S(\x_j,\x_k)\\ 
k>j&\si(i)<\si(k)<\si(j) && S(\x_k,\x_j) && S(\x_k,\x_i)
\end{array}\]

In all cases but the second the $S$-factors are exactly the same for $\s$ and $\s'$ when $\x_i=\x_j$. For the second we use $S(\x_k,\x)\,S(\x,\x_k)=1$. 

Sublemmas A.1--A.3 give Lemma A.1.

\section{Scaling Functions}

We give more detail into the derivation of the approximations (\ref{e:approx2}) and (\ref{e:approx3}), given in Section \ref{s:conjecture}. In the following, we fix the following two partitions
\begin{equation}
    J_1 \cup J_2 = K_1 \cup K_2 = \{1, 2, \dots, N \}
\end{equation}
so that $|J_i| = |K_i|$ for $i =1, 2$. Also, we fix a bijection $\tau: K_2 \rightarrow J_2$. Then, due to (\ref{e:approx1}), we will take the scaling
\begin{equation}\label{e:scaling2}
    x = -2t - st^{-1/3}, \quad \x_j = \mi + \mi\, \tilde{\x}_j\, t^{-1/3}, \quad \z_k = -\mi + \mi\, \tilde{\z}_k\, t^{-1/3}, \quad y_j +1= v_j \, t^{1/3} 
\end{equation}
for $j \in J_1$ and $k \in K_1$. Otherwise, we don't scale the variables.

\subsection{Approximation of $I_N(\x,\z;\tau)$}\label{s:I_approx}

We give more details for the approximation (\ref{e:approx2}). Take the definition of the function $I_N(\x, \z; \tau)$ given by (\ref{e:integrands}). In this case, since all the variables are labelled by the index sets $J_1$ and $K_1$, we take the scaling (\ref{e:scaling2}) for all the variables. The leading order of each term in $I_N(\x,\z;\tau)$ is given by the following:
\begin{equation}
    \begin{split}
    &\fr{\prod_{j \in J_1, k \in K_1}(\x_j + \z_k - 2\D \x_j \z_k)}{\prod_{\substack{j < k\\j,k \in J_1}}(1+ \x_j \x_k -2\D \x_j) \prod_{\substack{j<k\\j,k \in K_1}}(1 + \z_j \z_k - 2\D \z_j)} = (-2 \D)^{| J_1|} + \mathcal{O}(t^{-1/3})\\
    & d(\x_j, \z_{k}) = \fr{t^{1/3}}{(\z_{k} - \x_j)(- 2\D)} + \mathcal{O}(1)\\
    &D_N(\x, \z;\tau) = (-1)^{\sum_{k \in K_2}\tau(k) - k} \sum_{\gamma : K_1 \rightarrow J_1} (-2\D)^{-|J_1|}t^{|J_1| /3}\prod_{k \in K_1} \left( \fr{(-1)^{\gamma(k)- k}}{ \z_k -\x_{\gamma(k)}} + \mathcal{O}(t^{-1/3})\right)\\
    & \prod_{j \in J_1} \x_j^{x- y_j-1} e^{-\mi\, t \epsilon(\x_j) } = \prod_{j \in J_1}(\mi)^{x -y_j -1}\exp\left(\fr{1}{3} \tilde{\x}_j^3- (s +v_j)\,\tilde{\x}_j\,  \mathcal{O}(t^{-1/3}) \right)\\
    & \prod_{k \in K_1} \z_k^{x- y_k -1} e^{\mi \, t \, \epsilon(\z_k)} =\prod_{k \in K_1}(-\mi)^{x -y_k -1}\exp\left(-\fr{1}{3} \tilde{\z}_k^3+ (s+v_k)\,\tilde{\z}_j + \mathcal{O}(t^{-1/3}) \right).
    \end{split}
\end{equation}
Then, the approximation (\ref{e:approx2}) follows by combining the approximations above in the definition of $I_N(\x,\z;\tau)$ given by (\ref{e:integrands}).

\subsection{Approximation of $f(\x, \z;\tau)$}\label{s:f_approx}

We give more details for the approximation (\ref{e:approx3}). Take the definition of the function $f(\x, \z; \tau)$ given by (\ref{e:integrands}). In this case, all the $\z_j$ variable will be scaled as in (\ref{e:scaling2}), but only the $\x_k$ variables with $k\in J_1$ will be scaled as in (\ref{e:scaling2}). For simplicity, we take the following labelling of the index sets,
\begin{equation}
    K_2 = \{k_1 < k_2 < \cdots < k_M\}, \quad J_2 = \{\tau_1 = \tau(k_1), \dots, \tau_{M} = \tau(k_M) \}
\end{equation}
with $M = |K_2|$. Note that the last product term in $f(\x, \z;\tau)$ doesn't scale since all the variables are labeled by the index set $K_2$. Thus, we just consider the scaling of the other terms in $f(\x, \z;\tau)$.

Consider the term in $f(\x, \z; \tau)$,
\begin{equation}
    \prod_{k \in U(\ell)} \fr{1 + \x_{\tau_{\ell}} \x_{k} - 2 \D \x_k}{1 + \x_{\tau_{\ell}} \x_k - 2 \D \x_{\tau_{\ell}}}, \quad \text{with} \quad U(\ell) = \{k \mid  \tau_{\ell} < k, \quad k \neq \tau_{\ell+1}, \dots, \tau_{M}\}
\end{equation}
for a fixed $\ell$ with $1 \leq \ell \leq M$. The scaling of each term in the product will depend on the index $k$. If $k \in J_1$, then we take the scaling (\ref{e:scaling2}). Otherwise, if $k \in J_2$, then we don't take a scaling. Also, note that $k \in J_2$ if and only if $k = \tau_{\ell'}$ for some $\ell$ with $1 \leq \ell \leq M$. Then, we partition the set $U(\ell)=U_1(\ell) \cup U_{2}(\ell)$ with
\begin{equation}
    U_1(\ell) = \{k \mid  \tau_{\ell} < k, \quad k \neq \tau_{1}, \dots, \tau_{M}\}, \quad U_2(\ell) = \{k \mid  \tau_{\ell} < k, \quad k = \tau_{1}, \dots, \tau_{\ell-1}\}
\end{equation}
so that $U_1(\ell) \subset J_1$ and $U_2(\ell) \subset J_2$. Then, if $k\in U_1(\ell)$, we scale $\x_k$ as in (\ref{e:scaling2}) and obtain
\begin{equation}
    \fr{1 + \x_{\tau_{\ell}} \x_{k} - 2 \D \x_k}{1 + \x_{\tau_{\ell}} \x_k - 2 \D \x_{\tau_{\ell}}} = \fr{\mi\, \x_{\tau_{\ell}}- (\D \, \mi -1)}{(\mi - 2\D)\x_{\tau_{\ell}} +1} + \mathcal{O}(t^{-1/3}).
\end{equation}
On the other hand, if $k \in U_2(\ell)$, we don't scale the variables. Instead, we just note that $U_2(\ell) = \{\tau_{\ell'}  \mid \ell' < \ell, \quad \tau_{\ell} < \tau_{\ell'} \}$. Then, we have
\begin{equation}\label{e:approxU}
    \prod_{k \in U(\ell)} \fr{1 + \x_{\tau_{\ell}} \x_{k} - 2 \D \x_k}{1 + \x_{\tau_{\ell}} \x_k - 2 \D \x_{\tau_{\ell}}} = \left(\fr{\mi\, \x_{\tau_{\ell}}- (\D \, \mi -1)}{(\mi - 2\D)\x_{\tau_{\ell}} +1}  + \mathcal{O}(t^{-1/3})\right)^{|U_1(\ell)|} \prod_{k \in U_2(\ell)}\fr{1 + \x_{\tau_{\ell}} \x_{k} - 2 \D \x_k}{1 + \x_{\tau_{\ell}} \x_k - 2 \D \x_{\tau_{\ell}}}.
\end{equation}

Consider now the term in $f(\x, \z; \tau)$,
\begin{equation}
    \prod_{k \in V(\ell)} \fr{\x_{\tau_{\ell}} + \z_k - 2\D \x_{\tau_{\ell}} \z_{k}}{\x_{\tau_{\ell}} + \z_k - 2\D }, \quad \text{with} \quad V(\ell) = \{k \mid  k_{\ell} < k, \quad k \neq k_{\ell+1}, \dots, k_{M}\}
\end{equation}
for a fixed $\ell$ with $1 \leq \ell \leq M$. In this case, all the $\z_k$ variables will scale as in (\ref{e:scaling2}). Note that $V(\ell) = V_1(\ell)$ with
\begin{equation}
    V_1(\ell) = \{k \mid  k_{\ell} < k, \quad k \neq k_{1}, \dots, k_{M}\}
\end{equation}
since we have fixed the labelling so that $k_i < k_{i+1}$ for $i=1, \dots, M-1$. Then, when we scale the $\z_k$ variables, we obtain
\begin{equation} \label{e:approxV}
    \prod_{k \in V(\ell)} \fr{\x_{\tau_{\ell}} + \z_k - 2\D \x_{\tau_{\ell}} \z_{k}}{\x_{\tau_{\ell}} + \z_k - 2\D } = \left(\fr{(1 + 2\D\, \mi)\x_{\tau_{\ell}}- \mi}{\x_{\tau_{\ell}} - (2\D+\mi) } + \mathcal{O}(t^{-1/3}) \right)^{|V_1(\ell)|}.
\end{equation}

We now combine the approximations (\ref{e:approxU}) and (\ref{e:approxV}) and, then, take the product over $\ell$ to obtain the approximation (\ref{e:approx3}) for $f(\x, \z; \tau)$. Note that $|U_1(\ell)| - |V_1(\ell)| = \nu(k_{\ell}; \tau)$; see (\ref{e:nu1}). This follows by considering the following partitions,
\begin{equation}
    \begin{split}
    \{k \mid \tau_{\ell} <k\} &= U_1(\ell) \cup \{k \mid \tau_{\ell}<k, \quad k = \tau_1, \dots, \tau_M \}\\
    \{k \mid k_{\ell} <k\} &= V_1(\ell) \cup \{k \mid k_{\ell}<k, \quad k = k_1, \dots, k_M \}.
    \end{split}
\end{equation}
Additionally, note that we may write $U_2(\ell) = \{\tau_{\ell'}  \mid \ell' < \ell, \quad \tau_{\ell} < \tau_{\ell'} \}$. Then, we have
\begin{equation}
    \cup_{\ell=1}^M U_2(\ell) = \{(\tau_{\ell} < \tau_{\ell'}) \mid k_{\ell'} < k_{\ell} \}.
\end{equation}
Thus, the approximation (\ref{e:approx3}) follows by combining the approximations above in the definition of $f(\x,\z;\tau)$ given by (\ref{e:integrands}).


\end{document}